\documentclass[11pt,twoside]{article}
\usepackage{amsmath, amsthm, amscd, amsfonts, amssymb, graphicx, color}

\date{}

\begin{document}

\centerline{}

\centerline {\Large{\bf Continuous $K$-biframes in Hilbert spaces and their}}
\centerline {\Large{\bf  tensor products}}

\newcommand{\mvec}[1]{\mbox{\bfseries\itshape #1}}
\centerline{}
\centerline{\textbf{Prasenjit Ghosh}}
\centerline{Department of Mathematics, Barwan N. S. High School (HS),}
\centerline{Barwan, Murshidabad, 742161, West Bengal, India}
\centerline{e-mail: prasenjitpuremath@gmail.com}
\centerline{}
\centerline{\textbf{T. K. Samanta}}
\centerline{Department of Mathematics, Uluberia College,}
\centerline{Uluberia, Howrah, 711315,  West Bengal, India}
\centerline{e-mail: mumpu$_{-}$tapas5@yahoo.co.in}

\newtheorem{Theorem}{\quad Theorem}[section]

\newtheorem{definition}[Theorem]{\quad Definition}

\newtheorem{theorem}[Theorem]{\quad Theorem}

\newtheorem{remark}[Theorem]{\quad Remark}

\newtheorem{corollary}[Theorem]{\quad Corollary}

\newtheorem{note}[Theorem]{\quad Note}

\newtheorem{lemma}[Theorem]{\quad Lemma}

\newtheorem{example}[Theorem]{\quad Example}

\newtheorem{result}[Theorem]{\quad Result}
\newtheorem{conclusion}[Theorem]{\quad Conclusion}

\newtheorem{proposition}[Theorem]{\quad Proposition}

\begin{abstract}
\textbf{\emph{A generalization of continuous biframe in a Hilbert space is introduced and a few examples are discussed.\,Some characterizations and algebraic properties of this biframe are given.\,Here we also construct various types of continuous $K$-biframes with the help of a bounded linear operator.Relationship between continuous $K$-biframe and quotient operator is established.\,Finally, we define continuous $K$-biframe for the tensor products of Hilbert spaces and give an example.}}
\end{abstract}
\section{Introduction}
 
\smallskip\hspace{.6 cm}The notion of a frame in Hilbert space was first introduced by Duffin and Schaeffer \cite{Duffin} in connection with some fundamental problem in non-harmonic analysis.\,Thereafter, it was further developed and popularized by Daubechies et al.\,\cite{Daubechies} in 1986.\,A discrete frame is a countable family of elements in a separable Hilbert space which allows for a stable, not necessarily unique, decomposition of an arbitrary element into an expansion of the frame element.\,A sequence \,$\left\{\,f_{\,i}\,\right\}_{i \,=\, 1}^{\infty}$\, in a separable Hilbert space \,$H$\, is called a frame for \,$H$, if there exist positive constants \,$0 \,<\, A \,\leq\, B \,<\, \infty$\, such that
\[ A\; \|\,f\,\|^{\,2} \,\leq\, \sum\limits_{i \,=\, 1}^{\infty}\, \left|\ \left <\,f \,,\, f_{\,i} \, \right >\,\right|^{\,2} \,\leq\, B \,\|\,f\,\|^{\,2}\; \;\text{for all}\; \;f \,\in\, H.\]
The constants \,$A$\, and \,$B$\, are called lower and upper bounds, respectively.

Several generalizations of frames  namely, \,$g$-frame \cite{Sun}, \,$K$-frames \cite{Gavruta} etc. have been introduced in recent times.\;$K$-frames for a separable Hilbert spaces were introduced by Lara Gavruta to study the basic notions about atomic system for a bounded linear operator.\,Controlled frame \cite{I} is one of the newest generalization of frame in Hilbert space.\,I. Bogdanova et al.\,\cite{I} introduced controlled frame for spherical wavelets to get numerically more efficient approximation algorithm.\,Thereafter, P. Balaz \cite{B} developed weighted and controlled frame in Hilbert space.\,Biframe is also a generalization of controlled frame in Hilbert space which was studied by M. F. Parizi et al.\,\cite{MF}.\,To define frame in Hilbert space, only one sequence is needed, but for a biframe, two sequences are needed.\,A pair of sequences \,$\left(\,\left\{\,f_{i}\,\right\}_{i \,=\, 1}^{\,\infty}\,,\, \left\{\,g_{i}\,\right\}_{i \,=\, 1}^{\,\infty}\,\right)$\, in \,$H$\, is called a biframe for \,$H$\, if there exist positive constants \,$A$\, and \,$B$\, such that
\[A\; \|\,f\,\|^{\,2} \,\leq\, \sum\limits_{i \,=\, 1}^{\infty}\, \left <\,f\,,\, f_{\,i} \, \right >\,\left<\,g_{\,i}\,,\, f\,\right> \,\leq\, B \,\|\,f\,\|^{\,2}\; \;\forall\; f \,\in\, H.\] 
The constants \,$A$\, and \,$B$\, are called lower and upper biframe bounds, respectively.

\,Continuous frames extend the concept of discrete frames when the indices are related to some measurable space.\,Continuous frame in Hilbert space was studied by A. Rahimi et al \cite{AR}.\,M. H. faroughi and E. Osgooei \cite{MH} also studied continuous frame and continuous Bessel mapping.\,Continuous frame and discrete frame have been used in image processing, coding theory, wavelet analysis, signal denoising, feature extraction, robust signal processing etc.

In this paper, we give the notion of continuous $K$-biframe in Hilbert space, for some bounded linear operator \,$K$\, and then discuss some examples of this type of frame.\,A characterization of continuous $K$-biframe using its frame operator is established.\,We will see that the image of a continuous $K$-biframe under a bounded invertible operator in Hilbert space is also a continuous $K$-biframe in Hilbert space.\,Some algebraic properties of continuous $K$-biframe is established.\,We present a characterization of continuous $K$-biframe with the help of quotient operator.\,At the end, we discuss continuous $K$-biframes in tensor product of Hilbert spaces and  give a characterization. 

Throughout this paper,\;$H$\; is considered to be a Hilbert space with associated inner product \,$\left <\,\cdot,\, \cdot\,\right>$\, and \,$I_{H}$\; is the identity operator on \,$H$.\;$\mathcal{B}\,(\,H_{\,1},\, H_{\,2}\,)$\; is a collection of all bounded linear operators from \,$H_{\,1} \,\text{to}\, H_{\,2}$.\;In particular \,$\mathcal{B}\,(\,H\,)$\; denote the space of all bounded linear operators on \,$H$.\;For \,$T \,\in\, \mathcal{B}\,(\,H\,)$, we denote \,$\mathcal{N}\,(\,T\,)$\; and \,$\mathcal{R}\,(\,T\,)$\; for null space and range of \,$T$, respectively.\,$G\,\mathcal{B}\,(\,H\,)$\, denotes the set of all bounded linear operators which have bounded inverse.  
 
\section{Preliminaries}

\begin{theorem}\cite{O}\label{thm1.1}
Let \,$H_{\,1},\, H_{\,2}$\; be two Hilbert spaces and \;$U \,:\, H_{\,1} \,\to\, H_{\,2}$\; be a bounded linear operator with closed range \;$\mathcal{R}(\,U\,)$.\;Then there exists a bounded linear operator \,$U^{\dagger} \,:\, H_{\,2} \,\to\, H_{\,1}$\, such that \,$U\,U^{\dagger}\,x \,=\, x\; \;\forall\; x \,\in\, \mathcal{R}(\,U\,)$.
\end{theorem}

\begin{note}
The operator \,$U^{\dagger}$\, defined in Theorem (\ref{thm1.1}), is called the pseudo-inverse of \,$U$.
\end{note}

\begin{definition}{\cite{Ramu}}\label{def2.3}
Let \,$U,\, V \,\in\, \mathcal{B}\,(\,H\,)$\; with \,$\mathcal{N}\,(\,V\,) \,\subset\, \mathcal{N}\,(\,U\,)$.\;Then the linear operator \,$T \,=\, [\;U \,/\, V\;] \,:\, \mathcal{R}\,(\,V\,) \,\rightarrow\, \mathcal{R}\,(\,U\,)$, defined by \,$T\,(\,V\,f\,) \;=\; U\,f,\; f \,\in\, H$\, is called quotient operator on \,$H$. 
\end{definition}

\begin{note}
From definition \ref{def2.3}, it is easy to verify that \,$\mathcal{R}\,(\,T\,) \,\subset\, \mathcal{R}\,(\,U\,)$\, and \,$T\,V \,=\, U$. 
\end{note}

\begin{theorem}\cite{O}\label{2.5th1.001}
The set \,$\mathcal{S}\,(\,H\,)$\; of all self-adjoint operators on \,$H$\; is a partially ordered set with respect to the partial order \,$\leq$\, which is defined as for \,$T,\,S \,\in\, \mathcal{S}\,(\,H\,)$ 
\[T \,\leq\, S \,\Leftrightarrow\, \left<\,T\,f,\, f\,\right> \,\leq\, \left<\,S\,f,\, f\,\right>\; \;\forall\; f \,\in\, H.\] 
\end{theorem}

\begin{definition}{\cite{Gavruta}}
Let \,$K \,\in\, \mathcal{B}\,(\,H\,)$.\;Then a sequence \,$\{\,f_{\,i}\,\}_{i \,=\, 1}^{\infty}$\, in \,$H$\, is said to be a \,$K$-frame for \,$H$\, if there exist constants \,$0 \,<\, A \,\leq\, B \,<\, \infty$\, such that
\[A \,\left \|\,K^{\,\ast}\,f\, \right \|^{\,2} \,\leq\, \sum\limits^{\infty}_{i \,=\, 1}\, \left |\,\left <\,f \,,\, f_{\,i}\,\right >\,\right|^{\,2} \,\leq\, B\;\left\|\,f\,\right\|^{\,2}\; \;\forall\; f \,\in\, H.\]
\end{definition}

\begin{definition}\cite{AR}
Let \,$H$\, be a complex Hilbert space and \,$(\,\Omega,\, \mu\,)$\, be a measure space with positive measure \,$\mu$.\,A mapping \,$F \,:\, \Omega \,\to\, H$\, is called a continuous frame with respect to \,$\left(\,\Omega,\, \mu\,\right)$\, if
\begin{itemize}
\item[$(i)$] \,$F$\, is weakly-measurable, i.\,e., for all \,$f \,\in\, H$, \,$w \,\to\, \left<\,f,\, F\,(\,w\,)\,\right>$\, is a measurable function on \,$\Omega$,
\item[$(ii)$]there exist constants \,$0 \,<\, A \,\leq\, B \,<\, \infty$\, such that
\end{itemize}
\[A\,\left\|\,f\,\right\|^{\,2} \leq \int\limits_{\Omega}\left|\,\left<\,f,\, F\,(\,w\,)\,\right>\,\right|^{\,2}\,d\mu \leq B\left\|\,f\,\right\|^{\,2}\,,\]
for all \,$f \,\in\, H$.\,The constants \,$A$\, and \,$B$\, are called continuous frame bounds.\,If \,$A \,=\, B$, then it is called a tight continuous frame.\,If the mapping \,$F$\, satisfies only the right inequality, then it is called continuous Bessel mapping with Bessel bound \,$B$.
\end{definition}

Let \,$L^{\,2}\,(\,\Omega,\,\mu\,)$\, be the class of all measurable functions \,$f \,:\, \Omega \,\to\, H$\, such that \,$\|\,f\,\|^{\,2}_{\,2} \,=\,  \int\limits_{\,\Omega}\,\left\|\,f\,(\,w\,)\,\right\|^{\,2}\,d\mu \,<\, \infty$.\,It can be proved that \,$L^{\,2}\,(\,\Omega,\,\mu\,)$\, is a Hilbert space with respect to the inner product defined by
\[\left<\,f,\, g\right>_{L^{2}} \,=\, \int\limits_{\Omega}\,\left<\,f\,(\,w\,),\, g\,(\,w\,)\right>\,d\mu\,, \,f,\, g \,\in\, L^{\,2}\,(\,\Omega,\,\mu\,).\]

\begin{theorem}\cite{MH}
Let \,$F \,:\, \Omega \,\to\, H$\, be a Bessel mapping.\;Then the operator \,$T_{C} \,:\, L^{2}\left(\,\Omega,\,\mu\,\right) \,\to\, H$\, is defined by
\[\left<\,T_{C}\,(\,\varphi\,),\, h\,\right> \,=\, \int\limits_{\,\Omega}\,\varphi\,(\,w\,)\,\left<\,F\,(\,w\,),\, h\,\right>\,d\mu\]where \,$\varphi \,\in\, L^{2}\left(\,\Omega,\,\mu\,\right)$\, and \,$h \,\in\, H$\, is well-defined, linear, bounded and its adjoint operator is given by  
\[T^{\,\ast}_{C} \,:\, H \,\to\, L^{2}\left(\,\Omega,\,\mu\,\right) \;,\; T^{\,\ast}_{C}\,f\,(\,w\,) \,=\, \left<\,f,\, F\,(\,w\,)\,\right>\;,\; f \,\in\, H\,,\;\; w \,\in\, \Omega.\]
\end{theorem}

The operator \,$T_{C}$\, is called a pre-frame operator or synthesis operator and its adjoint operator is called analysis operator of \,$F$.\,

\begin{definition}\cite{MH}
Let \,$F \,:\, \Omega \,\to\, H$\, be a continuous frame for \,$H$.\,Then the operator \,$S_{C} \,:\, H \,\to\, H$\, defined by
\[\left<\,S_{C}\,(\,f\,),\, h\,\right> \,=\, \int\limits_{\,\Omega}\,\left<\,f, F\,(\,w\,)\,\right>\left<\,F\,(\,w\,),\, h\,\right>\,d\mu\,, \;\forall\, \,f,\, h \,\in\, H\]is called the frame operator of \,$F$.
\end{definition}

\begin{definition}\cite{Pro}
Let \,$H$\, be a Hilbert space and \,$(\,\Omega,\, \mu\,)$\, be a measure space with positive measure \,$\mu$.\,A pair \,$(\,\mathcal{F},\, \mathcal{G}\,) \,=\, \,\left(\,F \,:\, \Omega \,\to\, H,\, \,G \,:\, \Omega \,\to\, H\,\right)$\, of mappings is called a continuous biframe for \,$H$\, with respect to \,$\left(\,\Omega,\, \mu\,\right)$\, if
\begin{itemize}
\item[$(i)$] \,$F,\, G$\, are weakly-measurable, i.\,e., for all \,$f \,\in\, H$, \,$w \,\mapsto\, \left<\,f,\, F\,(\,w\,)\,\right>$\, and \,$w \,\mapsto\, \left<\,f,\, G\,(\,w\,)\,\right>$\, are measurable functions on \,$\Omega$,
\item[$(ii)$]there exist constants \,$0 \,<\, A \,\leq\, B \,<\, \infty$\, such that
\end{itemize}
\begin{align}
A\,\left\|\,f\,\right\|^{\,2} \leq \int\limits_{\Omega}\,\left<\,f,\, F\,(\,w\,)\,\right>\,\left<\,G\,(\,w\,),\, f\,\right>\,d\mu \leq\, B\left\|\,f\,\right\|^{\,2}\,,\label{3.eqq2.6}
\end{align}
for all \,$f \,\in\, H$.\,The constants \,$A$\, and \,$B$\, are called continuous biframe bounds.\,If \,$A \,=\, B$, then it is called a tight continuous biframe and it is called Parseval continuous biframe if \,$A \,=\, B \,=\, 1$.\,If the pair \,$(\,\mathcal{F},\, \mathcal{G}\,)$\, satisfies only the right inequality, then it is called continuous biframe Bessel mapping with Bessel bound \,$B$. 
\end{definition}

In particular, if \,$\mu$\, is a counting measure and \,$\Omega \,=\, \mathbb{N}$, then \,$(\,\mathcal{F},\, \mathcal{G}\,)$\, is called a discrete biframe for \,$H$.

\begin{definition}\cite{Pro}
Let \,$(\,\mathcal{F},\, \mathcal{G}\,) \,=\, \,\left(\,F \,:\, \Omega \,\to\, H,\, \,G \,:\, \Omega \,\to\, H\,\right)$\, be a continuous biframe for \,$H$\, with respect to \,$\left(\,\Omega,\, \mu\,\right)$.\,Then the operator \,$S_{\mathcal{F},\, \mathcal{G}} \,:\, H \,\to\, H$\, defined by
\[S_{\mathcal{F},\, \mathcal{G}}\,f \,=\, \int\limits_{\,\Omega}\,\left<\,f, F\,(\,w\,)\,\right>\,G\,(\,w\,)\,d\mu,\]
for all \,$f \,\in\, H$\, is called the frame operator.
\end{definition}

Now, for each \,$f \,\in\, H$, we have 
\begin{align}
\left<\,S_{\mathcal{F},\, \mathcal{G}}\,f,\, f\,\right> \,=\, \int\limits_{\,\Omega}\,\left<\,f, F\,(\,w\,)\,\right>\left<\,G\,(\,w\,),\, f\,\right>\,d\mu.\label{3.eqq3.12}
\end{align}
Thus, for each \,$f \,\in\, H$, we get
\[A\,\left\|\,f\,\right\|^{\,2} \leq\, \left<\,S_{\mathcal{F},\, \mathcal{G}}\,f,\, f\,\right> \,\leq\,  B\left\|\,f\,\right\|^{\,2}.\]
This implies that \,$A\,I_{H} \,\leq\, S_{\mathcal{F},\, \mathcal{G}} \,\leq\, B\,I_{H}$, where \,$I_{H}$\, is the identity operator on \,$H$.\,Hence, \,$S_{\mathcal{F},\, \mathcal{G}}$\, is positive and invertible.\,Here, we assume that  \,$S_{\mathcal{F},\, \mathcal{G}}$\, is self-adjoint operator.\\

The tensor product of Hilbert spaces are introduced by several ways and it is a certain linear space of operators which was represented by Folland in \cite{Folland}.

\begin{definition}\cite{Upender}\label{1.def1.01}
The tensor product of Hilbert spaces \,$\left(\,H_{1},\, \left<\,\cdot,\,\cdot\,\right>_{1}\,\right)$\, and \,$\left(\,H_{2},\, \left<\,\cdot,\,\cdot\,\right>_{2}\,\right)$\, is denoted by \,$H_{1} \,\otimes\, H_{2}$\, and it is defined to be an inner product space associated with the inner product
\begin{equation}\label{eq1.001}   
\left<\,f \,\otimes\, g \,,\, f^{\,\prime} \,\otimes\, g^{\,\prime}\,\right> \,=\, \left<\,f,\, f^{\,\prime}\,\right>_{\,1}\;\left<\,g,\, g^{\,\prime}\,\right>_{\,2},
\end{equation}
for all \,$f,\, f^{\,\prime} \,\in\, H_{1}\; \;\text{and}\; \;g,\, g^{\,\prime} \,\in\, H_{2}$.\,The norm on \,$H_{1} \,\otimes\, H_{2}$\, is given by 
\begin{equation}\label{eq1.0001}
\left\|\,f \,\otimes\, g\,\right\| \,=\, \|\,f\,\|_{\,1}\;\|\,g\,\|_{\,2}\; \;\forall\; f \,\in\, H_{1}\; \;\text{and}\; \,g \,\in\, H_{2}.
\end{equation}
The space \,$H_{1} \,\otimes\, H_{2}$\, is complete with respect to the above inner product.\;Therefore the space \,$H_{1} \,\otimes\, H_{2}$\, is a Hilbert space.     
\end{definition} 

For \,$Q \,\in\, \mathcal{B}\,(\,H_{1}\,)$\, and \,$T \,\in\, \mathcal{B}\,(\,H_{2}\,)$, the tensor product of operators \,$Q$\, and \,$T$\, is denoted by \,$Q \,\otimes\, T$\, and defined as 
\[\left(\,Q \,\otimes\, T\,\right)\,A \,=\, Q\,A\,T^{\,\ast}\; \;\forall\; \;A \,\in\, H_{1} \,\otimes\, H_{2}.\]
It can be easily verified that \,$Q \,\otimes\, T \,\in\, \mathcal{B}\,(\,H_{1} \,\otimes\, H_{2}\,)$\, \cite{Folland}.\\

\begin{theorem}\cite{Folland}\label{th1.1}
Suppose \,$Q,\, Q^{\prime} \,\in\, \mathcal{B}\,(\,H_{1}\,)$\, and \,$T,\, T^{\prime} \,\in\, \mathcal{B}\,(\,H_{2}\,)$, then \begin{itemize}
\item[$(i)$] \,$Q \,\otimes\, T \,\in\, \mathcal{B}\,(\,H_{1} \,\otimes\, H_{2}\,)$\, and \,$\left\|\,Q \,\otimes\, T\,\right\| \,=\, \|\,Q\,\|\; \|\,T\,\|$.
\item[$(ii)$] \,$\left(\,Q \,\otimes\, T\,\right)\,(\,f \,\otimes\, g\,) \,=\, Q\,(\,f\,) \,\otimes\, T\,(\,g\,)$\, for all \,$f \,\in\, H_{1},\, g \,\in\, H_{2}$.
\item[$(iii)$] $\left(\,Q \,\otimes\, T\,\right)\,\left(\,Q^{\,\prime} \,\otimes\, T^{\,\prime}\,\right) \,=\, (\,Q\,Q^{\,\prime}\,) \,\otimes\, (\,T\,T^{\,\prime}\,)$. 
\item[$(iv)$] \,$Q \,\otimes\, T$\, is invertible if and only if \,$Q$\, and \,$T$\, are invertible, in which case \,$\left(\,Q \,\otimes\, T\,\right)^{\,-\, 1} \,=\, \left(\,Q^{\,-\, 1} \,\otimes\, T^{\,-\, 1}\,\right)$.
\item[$(v)$] \,$\left(\,Q \,\otimes\, T\,\right)^{\,\ast} \,=\, \left(\,Q^{\,\ast} \,\otimes\, T^{\,\ast}\,\right)$.  
\end{itemize}
\end{theorem}

\section{Continuous $K$-biframe in Hilbert space}
In this section, first we give the definition of a continuous $K$-biframe in Hilbert space and then discuss some properties.

\begin{definition}\label{def1.01}
Let \,$H$\, be a Hilbert space and \,$K$\, be a bounded linear operator on \,$H$.\,Suppose \,$(\,\Omega,\, \mu\,)$\, is a measure space with positive measure \,$\mu$.\,A pair \,$(\,\mathcal{F},\, \mathcal{G}\,) \,=\, \,\left(\,F \,:\, \Omega \,\to\, H,\, \,G \,:\, \Omega \,\to\, H\,\right)$\, of mappings is called a continuous $K$-biframe for \,$H$\, with respect to \,$\left(\,\Omega,\, \mu\,\right)$\, if
\begin{itemize}
\item[$(i)$] \,$F,\, G$\, are weakly-measurable, i.\,e., for all \,$f \,\in\, H$, \,$w \,\mapsto\, \left<\,f,\, F\,(\,w\,)\,\right>$\, and \,$w \,\mapsto\, \left<\,G\,(\,w\,),\, f\,\right>$\, are measurable functions on \,$\Omega$,
\item[$(ii)$]there exist constants \,$0 \,<\, A \,\leq\, B \,<\, \infty$\, such that
\end{itemize}
\[A\,\left\|\,K^{\,\ast}\,f\,\right\|^{\,2} \leq \int\limits_{\Omega}\,\left<\,f,\, F\,(\,w\,)\,\right>\,\left<\,G\,(\,w\,),\, f\,\right>\,d\mu \leq\, B\left\|\,f\,\right\|^{\,2}\,,\]
for all \,$f \,\in\, H$.\,The constants \,$A$\, and \,$B$\, are called continuous $K$-biframe bounds.\,If \,$A \,=\, B$, then it is called a tight continuous $K$-biframe and it is called Parseval continuous $K$-biframe if \,$A \,=\, B \,=\, 1$.\,If the pair \,$(\,\mathcal{F},\, \mathcal{G}\,)$\, satisfies only the right inequality, then it is called continuous biframe Bessel mapping with Bessel bound \,$B$.\,If \,$K \,=\, I_{H}$\, and \,$\mu$\, is a counting measure with \,$\Omega \,=\, \mathbb{N}$, then \,$(\,\mathcal{F},\, \mathcal{G}\,)$\, is called a discrete biframe for \,$H$.   
\end{definition}

\begin{remark}
Let \,$F \,:\, \Omega \,\to\, H$\, be a mapping.\,Then according to the definition \ref{def1.01}, we say that
\begin{itemize}
\item[$(i)$]If \,$(\,\mathcal{F},\, \mathcal{F}\,)$\, is a continuous $K$-biframe for \,$H$, then \,$\mathcal{F}$\, is a continuous $K$-frame for \,$H$.
\item[$(ii)$]If for some \,$U \,\in\, G\,\mathcal{B}(\,H\,)$, \,$(\,\mathcal{F},\, U\,\mathcal{F}\,)$\, is a continuous $K$-biframe for \,$H$, then \,$\mathcal{F}$\, is a \,$U$-controlled continuous $K$-frame for \,$H$.  
\item[$(iii)$]If for some \,$T,\, U \,\in\, G\,\mathcal{B}(\,H\,)$, \,$(\,T\,\mathcal{F},\, U\,\mathcal{F}\,)$\, is a continuous $K$-biframe for \,$H$, then \,$\mathcal{F}$\, is a \,$(\,T,\, U\,)$-controlled continuous $K$-frame for \,$H$.  
\end{itemize} 
\end{remark}

Now, we validate the above definition by some examples.

\begin{example}
Let \,$H \,=\, \mathbb{R}^{\,3}$\, and \,$\left\{\,e_{\,1},\,e_{\,2},\, e_{\,3}\,\right\}$\, be an standard orthonormal basis for \,$H$.\,Consider \,$\Omega \,=\, \left\{\,x \,\in\, \mathbb{R}^{\,3} \,:\, \|\,x\,\| \,\leq\, 1\,\right\}$.\,Then it is a measure space equipped with the Lebesgue measure \,$\mu$.\,Suppose \,$\left\{\,B_{\,1},\,B_{\,2},\, B_{\,3}\,\right\}$\, is a partition of \,$\Omega$\, where \,$\mu\,(\,B_{1}\,) \,\geq\, \mu\,(\,B_{2}\,) \,\geq\, \mu\,(\,B_{3}\,) \,>\, 1$.\,Define 
\[F \,:\, \Omega \,\to\, H\hspace{.5cm}\text{by} \hspace{.3cm} F\,(\,w\,) \,=\, \begin{cases}
\dfrac{2\,e_{\,1}}{\sqrt{\,\mu\left(\,B_{1}\,\right)}} & \text{if\;\;}\; w \,\in\, B_{1} \\ \dfrac{3\,e_{\,2}}{\sqrt{\,\mu\left(\,B_{2}\,\right)}} & \text{if\;\;}\; w \,\in\, B_{2}\\ \dfrac{\,-\,2\,e_{\,3}}{\sqrt{\,\mu\left(\,B_{3}\,\right)}} & \text{if\;\;}\; w \,\in\, B_{3}\,, \end{cases}\]
\[G \,:\, \Omega \,\to\, H\hspace{.5cm}\text{by} \hspace{.3cm} G\,(\,w\,) \,=\, \begin{cases}
\dfrac{2\,e_{\,1}}{\sqrt{\,\mu\left(\,B_{1}\,\right)}} & \text{if\;\;}\; w \,\in\, B_{1} \\ \dfrac{\, e_{\,2}}{\sqrt{\,\mu\left(\,B_{2}\,\right)}} & \text{if\;\;}\; w \,\in\, B_{2}\\ \dfrac{-\,e_{\,3}}{\sqrt{\,\mu\left(\,B_{3}\,\right)}} & \text{if\;\;}\; w \,\in\, B_{3} \end{cases}\]
It is easy to verify that for all \,$f \,\in\, H$, \,$w \,\mapsto\, \left<\,f,\, F\,(\,w\,)\,\right>$\, and \,$w \,\mapsto\, \left<\,f,\, G\,(\,w\,)\,\right>$\, are measurable functions on \,$\Omega$.\,Define \,$K \,:\, H \,\to\, H$\, by \,$K\,e_{1} \,=\, e_{1}$,\, \,$K\,e_{2} \,=\, e_{3}$\, and \,$K\,e_{3} \,=\, e_{2}$.\,Then \,$K^{\,\ast}\,e_{1} \,=\, e_{1}$,\, \,$K^{\,\ast}\,e_{2} \,=\, e_{3}$, \,$K^{\,\ast}\,e_{3} \,=\, e_{2}$\, and it is easy to verify that \,$\left\|\,K^{\,\ast}\,f\,\right\|^{\,2} \,=\, \left\|\,f\,\right\|^{\,2}$.\,Now, for \,$f \,\in\, H$, we have
\begin{align*}
\int\limits_{\,\Omega}\left|\left<\,f,\, F\,(\,w\,)\,\right>\,\right|^{2}\,d\mu &= \int\limits_{\,B_{1}}\left|\,\left<\,f,\, \dfrac{2}{\sqrt{\mu\left(\,B_{1}\,\right)}}\,e_{\,1}\,\right>\,\right|^{\,2}\,d\mu \,+\, \int\limits_{\,B_{2}}\left|\,\left<\,f,\, \dfrac{3}{\sqrt{\mu\left(\,B_{2}\,\right)}}\,e_{\,2}\,\right>\,\right|^{\,2}\,d\mu\\
&\hspace{1cm}+\,\int\limits_{\,B_{3}}\left|\,\left<\,f,\, \dfrac{-\,1}{\sqrt{\mu\left(\,B_{3}\,\right)}}\,2\,e_{\,1}\,\right>\,\right|^{\,2}\,d\mu\\
&=\,4\,\left|\,\left<\,f,\, \,e_{\,1}\,\right>\,\right|^{\,2} \,+\, 9\,\left|\,\left<\,f,\, \,e_{\,2}\,\right>\,\right|^{\,2} \,+\, 4\,\left|\,\left<\,f,\, \,e_{\,3}\,\right>\,\right|^{\,2}\\
& \,=\, 4\,\|\,f\,\|^{\,2} \,+\, 5\,\left|\,\left<\,f,\, \,e_{\,3}\,\right>\,\right|^{\,2}.
\end{align*} 
Therefore, \,$F \,:\, \Omega \,\to\, H$\, is a continuous $K$-frame for \,$H$\, with bounds \,$4$\, and \,$9$.\,Similarly, it can be shown that \,$G \,:\, \Omega \,\to\, H$\, is a continuous $K$-frame for \,$H$\, with bounds \,$1$\, and \,$4$.\,On the other hand, for \,$f \,\in\, H$, we have
\begin{align*}
&\int\limits_{\,\Omega}\,\left<\,f,\, F\,(\,w\,)\,\right>\,\left<\,G\,(\,w\,),\, f\,\right>\,d\mu \\
&= \int\limits_{\,B_{1}}\left<\,f,\, \dfrac{2\,e_{\,1}}{\sqrt{\mu\left(\,B_{1}\,\right)}}\,\right>\,\left<\,\dfrac{2\,e_{\,1}}{\sqrt{\mu\left(\,B_{1}\,\right)}},\, f\,\right>\,d\mu\\
& \,+\, \int\limits_{\,B_{2}}\left<\,f,\, \dfrac{3\,e_{\,2}\,}{\sqrt{\mu\left(\,B_{2}\,\right)}}\,\right>\,\left<\,\dfrac{\, e_{\,2}}{\sqrt{\mu\left(\,B_{2}\,\right)}},\, f\,\right>\,d\mu\\
&+\,\int\limits_{\,B_{3}}\left<\,f,\, \dfrac{-\,2\,e_{\,3}\,}{\sqrt{\mu\left(\,B_{3}\,\right)}}\,\right>\,\left<\, \dfrac{-\,e_{\,3}\,}{\sqrt{\mu\left(\,B_{3}\,\right)}},\, f\,\right>\,d\mu\\
&=\,4\,\left|\,\left<\,f,\, \,e_{\,1}\,\right>\,\right|^{\,2} \,+\, 3\,\left|\,\left<\,f,\, \,e_{\,2}\,\right>\,\right|^{\,2} \,+\, 2\,\left|\,\left<\,f,\, \,e_{\,3}\,\right>\,\right|^{\,2}\\
& \,=\, 2\,\|\,f\,\|^{\,2} \,+\, 2\,\left|\,\left<\,f,\, \,e_{\,1}\,\right>\,\right|^{\,2} \,+\, \left|\,\left<\,f,\, \,e_{\,2}\,\right>\,\right|^{\,2}.
\end{align*} 
Thus, for each \,$f \,\in\, H$, we get
\begin{align*}
2\,\left\|\,K^{\,\ast}\,f\,\right\|^{\,2} \,\leq\, \int\limits_{\,\Omega}\,\left<\,f,\, F\,(\,w\,)\,\right>\,\left<\,G\,(\,w\,),\, f\,\right>\,d\mu \,\leq\, 5\,\|\,f\,\|^{\,2}. 
\end{align*}
Therefore, \,$(\,\mathcal{F},\, \mathcal{G}\,)$\, is a continuous $K$-biframe for \,$H$\, with bounds \,$2$\, and \,$5$. 
\end{example}

\begin{example}
Let \,$H \,=\, \mathbb{R}^{\,3}$\, and \,$\left\{\,e_{\,1},\,e_{\,2},\, e_{\,3}\,\right\}$\, be an standard orthonormal basis for \,$H$.\,Define \,$K \,:\, \mathbb{R}^{\,3} \,\to\, \mathbb{R}^{\,3}$\, by \,$K\left(\,x,\, y,\, z\,\right) \,=\, \left(\,x \,+\, y,\,  y \,+\, z,\, z \,+\, x\,\right)$, for every \,$\left(\,x,\, y,\, z\,\right) \,\in\, \mathbb{R}^{\,3}$.\,Then \,$K$\, is bounded linear operator on \,$\mathbb{R}^{\,3}$.\,The matrix associated with the operator \,$K$\, with respect to the above basis is given by
\[
\left[\,K\,\right] \,=\, 
\begin{pmatrix}
\;1 & 0 & 1\\
\; 1   & 1 & \,0\\
\; 0    & 1  & 1
\end{pmatrix}
.\] 
Thus, for \,$\left(\,x,\, y,\, z\,\right) \,\in\, \mathbb{R}^{\,3}$, we get \,$K^{\,\ast}\left(\,x,\, y,\, z\,\right) \,=\, \left(\,x \,+\, z,\,  x \,+\, y,\, y \,+\, z\,\right)$\, and 
\begin{align*}
\left\|\,K^{\,\ast}\,\left(\,x,\, y,\, z\,\right)\,\right\|^{\,2} &\,=\, \left(\,x \,+\, z\,\right)^{2} \,+\, \left(\,x \,+\, y\,\right)^{2} \,+\, \left(\,y \,+\, z\,\right)^{2} \\
&=\, 2\,\left(\,x^{\,2} \,+\, y^{\,2} \,+\, z^{\,2}\,\right) \,+\, 2\,\left(\,x\,y \,+\, y\,z \,+\, z\,x\,\right)\\
&\leq\,4\,\left(\,x^{\,2} \,+\, y^{\,2} \,+\, z^{\,2}\,\right).
\end{align*}
Now, we consider a measure space \,$\left(\,\Omega \,=\, [\,0,\,1\,],\, \mu\,\right)$\, where \,$\mu$\,  is the Lebesgue measure.\,For each \,$w \,\in\, \Omega$, we define the mappings \,$F \,:\, \Omega \,\to\, H$\, by \,$F\,(\,w\,) \,=\, \left(\,2\,w,\, 1,\, 1\,\right)$\, and \,$G \,:\, \Omega \,\to\, H$\, by \,$G\,(\,w\,) \,=\, \left(\,1,\, 2\,w,\, 1\,\right)$.\,Then for \,$f \,=\, \left(\,x,\, y,\, z\,\right) \,\in\, H$, we have
\begin{align*}
&\int\limits_{\,\Omega}\,\left<\,f,\, F\,(\,w\,)\,\right>\,\left<\,G\,(\,w\,),\, f\,\right>\,d\mu \\
&=\,\int\limits_{[\,0,\, 1\,]}\,\left<\,\left(\,x,\, y,\, z\,\right),\, \left(\,2\,w,\, 1,\, 1\,\right)\,\right>\,\left<\,\left(\,1,\, 2\,w,\, 1\,\right),\, \left(\,x,\, y,\, z\,\right)\,\right>\,d\mu \\
&=\,\int\limits_{[\,0,\, 1\,]}\,\left(\,2\,x\,w \,+\, y \,+\, z\,\right)\,\left(\,x \,+\, 2\,y\,w \,+\, z\,\right)\,d\mu \\
&=\,x^{\,2} \,+\, y^{\,2} \,+\, z^{\,2} \,+\, \dfrac{7}{3}\,x\,y \,+\, 2\,x\,z \,+\, 2\,y\,z\\
&\leq\,3\left(\,x^{\,2} \,+\, y^{\,2} \,+\, z^{\,2}\,\right) \,+\, \dfrac{x\,y}{3} \,\leq\,3\left(\,x^{\,2} \,+\, y^{\,2} \,+\, z^{\,2}\,\right) \,+\, \dfrac{x^{\,2} \,+\, y^{\,2}}{6}.     
\end{align*}
Thus, for each \,$f \,\in\, H$, we get
\begin{align*}
&3\,\,\left\|\,f\,\right\|^{\,2} \leq \int\limits_{\Omega}\,\left<\,f,\, F\,(\,w\,)\,\right>\,\left<\,G\,(\,w\,),\, f\,\right>\,d\mu \leq\, 4\,\left\|\,f\,\right\|^{\,2}.\\
&\Rightarrow\,\dfrac{3}{4}\left\|\,K^{\,\ast}\,f\,\right\|^{\,2} \leq \int\limits_{\Omega}\,\left<\,f,\, F\,(\,w\,)\,\right>\,\left<\,G\,(\,w\,),\, f\,\right>\,d\mu \leq\, 4\left\|\,f\,\right\|^{\,2}.
\end{align*}
\,$(\,\mathcal{F},\, \mathcal{G}\,)$\, is a continuous $K$-biframe for \,$H$\, with bounds \,$3 \,/\, 4$\, and \,$4$. 
\end{example}

\begin{example}\label{3exm3.11}
Let \,$H$\, be an infinite dimensional separable Hilbert space and \,$\left\{\,e_{\,i}\,\right\}_{i \,=\, 1}^{\,\infty}$\, be an orthonormal basis for \,$H$.\,Suppose
\begin{align*}
&\left\{\,f_{\,i}\,\right\}_{i \,=\, 1}^{\,\infty} \,=\, \left\{\,3\,e_{\,1},\, 2\,e_{\,2},\, 2\,e_{\,3},\, \cdots\,\cdots\,\right\}\,,\\
&\left\{\,g_{\,i}\,\right\}_{i \,=\, 1}^{\,\infty} \,=\, \left\{\,0,\, 2\,e_{\,1},\, 0,\, e_{\,2},\, 0,\, e_{\,3},\, \cdots\,\cdots\,\right\}.
\end{align*}
Define \,$K \,:\, H \,\to\, H$\, by \,$K\,f \,=\, \sum\limits_{i \,=\, 1}^{\,m}\,\left<\,f,\, e_{\,i}\,\right>\,e_{\,i}$, \,$f \,\in\, H$, where \,$m$\, is a fixed positive integer.\,Now, for each \,$f \,\in\, H$, we have
\[\left \|\,K^{\,\ast}\,f\, \right \|^{\,2} \,=\, \sum\limits_{i \,=\, 1}^{\,m}\,\left|\,\left<\,f,\, e_{\.i}\,\right>\,\right|^{\,2} \,\leq\, \sum\limits_{i \,=\, 1}^{\,\infty}\,\left|\,\left<\,f,\, e_{\.i}\,\right>\,\right|^{\,2} \,=\, \left \|\,f\, \right \|^{\,2}.\]
Now, let \,$\left(\,\Omega,\, \mu\,\right)$\, be a measure space with \,$\mu$\, is \,$\sigma$-finite.\,Then we can write \,$\Omega \,=\, \bigcup_{i \,=\, 1}^{\,\infty}\,\Omega_{i}$, where \,$\left\{\,\Omega_{\,i}\,\right\}_{i \,=\, 1}^{\,\infty}$\, is a sequence of disjoint measurable subsets of \,$\Omega$\, with \,$\mu\left(\,\Omega_{i}\,\right) \,<\, \infty$.\,For each \,$w \,\in\, \Omega$, we define the mappings \,$F \,:\, \Omega \,\to\, H$\, by \,$F\,(\,w\,) \,=\, \dfrac{1}{\sqrt{\mu\left(\,\Omega_{i}\,\right)}}\,f_{\,i}$\, and \,$G \,:\, \Omega \,\to\, H$\, by \,$G\,(\,w\,) \,=\, \dfrac{1}{\sqrt{\mu\left(\,\Omega_{i}\,\right)}}\,g_{\,i}$.\,Then for \,$f \,\in\, H$, we have
\begin{align*}
\int\limits_{\,\Omega}\left|\left<\,f,\, F\,(\,w\,)\,\right>\,\right|^{2}\,d\mu& \,=\,\sum\limits_{i \,=\, 1}^{\infty}\int\limits_{\,\Omega_{i}}\left|\,\left <\,f,\, f_{\,i}\,\right >\,\right|^{\,2}\,d\mu \,=\, \left|\,\left <\,f,\, e_{\,1}\,\right >\,\right|^{\,2} \,+\, 2\,\sum\limits_{i \,=\, 1}^{\infty}\, \left|\,\left <\,f,\, e_{\,i} \,\right >\,\right|^{\,2}\\
&=\,\left|\,\left <\,f,\, e_{\,1}\,\right >\,\right|^{\,2} \,+\, 2\,\left\|\,f\,\right\|^{\,2}.
\end{align*} 
Therefore, \,$F$\, is a continuous $K$-frame for \,$H$\, with bounds \,$2$\, and \,$3$.\,Similarly, it can be shown that \,$G$\, is a continuous $K$-frame for \,$H$\, with bounds \,$1$\, and \,$2$.\,Now, for \,$f \,\in\, H$, we have
\begin{align*}
&\int\limits_{\,\Omega}\,\left<\,f,\, F\,(\,w\,)\,\right>\,\left<\,G\,(\,w\,),\, f\,\right>\,d\mu \,=\,\sum\limits_{i \,=\, 1}^{\infty}\int\limits_{\,\Omega_{i}}\,\left <\,f,\, f_{\,i}\,\right >\,\left <\,g_{\,i},\, f\,\right >\,d\mu\\
&=\,\left <\,f,\, e_{\,1}\,\right >\,\left <\,e_{\,1},\, f\,\right > \,+\, \left <\,f,\, e_{\,1}\,\right >\,\left <\,e_{\,1},\, f\,\right > \,+\, \left <\,f,\, e_{\,2}\,\right >\,\left <\,e_{\,2},\, f\,\right >+\cdots\\
&=\,\left|\,\left <\,f,\, e_{\,1}\,\right >\,\right|^{\,2} \,+\, \left|\,\left <\,f,\, e_{\,1}\,\right >\,\right|^{\,2} \,+\, \left|\,\left <\,f,\, e_{\,2}\,\right >\,\right|^{\,2} \,+\, \left|\,\left <\,f,\, e_{\,3}\,\right >\,\right|^{\,2}\,+\,\cdots\\
& \,=\, \left|\,\left <\,f,\, e_{\,1}\,\right >\,\right|^{\,2} \,+\, \|\,f\,\|^{\,2}. 
\end{align*} 
Thus, for each \,$f \,\in\, H$, we get
\begin{align*}
\left\|\,K^{\,\ast}\,f\,\right\|^{\,2} \,\leq\, \int\limits_{\,\Omega}\,\left<\,f,\, F\,(\,w\,)\,\right>\,\left<\,G\,(\,w\,),\, f\,\right>\,d\mu \,\leq\, 2\,\|\,f\,\|^{\,2}. 
\end{align*}
Hence, \,$(\,\mathcal{F},\, \mathcal{G}\,)$\, is a continuous $K$-biframe for \,$H$\, with bounds \,$1$\, and \,$2$.
\end{example}

Next, we give an example of a continuous $K$-biframe for a real inner product space.

\begin{example}
Let
\[
V \,=\, \left\{
\begin{pmatrix}
\;a & 0\\
\;0           & b\\
\end{pmatrix}
:\, a,\, b \,\in\, \mathbb{R}
\right\}.
\]
Define \,$\left<\,\cdot,\, \cdot\,\right> \,:\, V \,\times\, V \,\to\, \mathbb{R}$\, by \,$\left<\,M,\, N\,\right> \,=\, det\left(\,MN^{\,t}\,\right)$, for all \,$M,\, N \,\in V$.\,Then it is easy to verify that \,$\left<\,\cdot,\, \cdot\,\right>$\, is an real inner product on \,$V$.\,Now, we consider a measure space \,$\left(\,\Omega \,=\, [\,0,\,1\,],\, \mu\,\right)$\, where \,$\mu$\,  is the Lebesgue measure.\,Define \,$F \,:\, \Omega \,\to\, V$\, by
\[
F\,(\,w\,) \,=\, 
\begin{pmatrix}
\;(\,\sqrt{3} \,+\, 1\,)\,(\,1 \,-\, w\,) & 0\\
\;0           & (\,\sqrt{3} \,-\, 1\,)\,(\,1 \,+\, w\,)\\
\,\end{pmatrix}
,\; w \,\in\, \Omega
\] 
and \,$G \,:\, \Omega \,\to\, V$\, by
\[
G\,(\,w\,) \,=\, 
\begin{pmatrix}
\;\sqrt{11\,w}\, & 0\\
\;0           & \sqrt{11\,w}\,\\
\end{pmatrix}
,\; w \,\in\, \Omega
\]
It is easy to verify that for all \,$M \,\in\, V$, \,$w \,\mapsto\, \left<\,M,\, F\,(\,w\,)\,\right>$\, and \,$w \,\mapsto\, \left<\,M,\, G\,(\,w\,)\,\right>$\, are measurable functions on \,$\Omega$.\,Define \,$K \,:\, V \,\to\, V$\, by \,$K\,M \,=\, 2\,M$,\, for all \,$M \,\in\, V$.\,Then it is easy to verify that \,$\left\|\,K^{\,\ast}\,M\,\right\|^{\,2} \,=\, 4\,(\,det\,M\,)^{\,2} \,=\, 4\,\left\|\,M\,\right\|^{\,2}$.\,Now, for each \,$
M \,=\, 
\begin{pmatrix}
\;a\, & 0\\
\;0           & b\,\\
\end{pmatrix}
 \,\in\, V
$, we get
\begin{align*}
\left<\,M,\, F\,(\,w\,)\,\right> \,=\, det\left\{\,\begin{pmatrix}
\;a\, & 0\\
\;0           & b\,\\
\end{pmatrix}
\begin{pmatrix}
\;(\,\sqrt{3} \,+\, 1\,)\,(\,1 \,-\, w\,) & 0\\
\;0           & (\,\sqrt{3} \,-\, 1\,)\,(\,1 \,+\, w\,)\\
\end{pmatrix}
\,\right\} \\
\,=\, 2\,a\,b\,\left(\,1 \,-\, w^{\,2}\,\right)\hspace{7.9cm}
\end{align*} 
and
\begin{align*}
\left<\,G\,(\,w\,),\, M\,\right> \,=\, det\left\{\,\begin{pmatrix}
\;\sqrt{11\,w}\, & 0\\
\;0           & \sqrt{11\,w}\,\\
\end{pmatrix}
\begin{pmatrix}
\;a\, & 0\\
\;0           & b\,\\
\end{pmatrix}
\,\right\} \,=\, 11\,a\,b\,w.
\end{align*} 
Thus, for each \,$
M \,=\, 
\begin{pmatrix}
\;a\, & 0\\
\;0           & b\,\\
\end{pmatrix}
 \,\in\, V
$, we have
\begin{align*}
&\int\limits_{\,[\,0,\, 1\,]}\,\left<\,M,\, F\,(\,w\,)\,\right>\,\left<\,G\,(\,w\,),\, M\,\right>\,d\mu \,=\, \int\limits_{\,[\,0,\, 1\,]}\,22\,a^{\,2}\,b^{\,2}\,w\,\left(\,1 \,-\, w^{\,2}\,\right)\,d\mu \\
&\hspace{2.5cm}=\,\dfrac{11}{2}\,a^{\,2}\,b^{\,2} \,=\, \dfrac{11}{2}\,\,det \begin{pmatrix}
\;a^{\,2}\, & 0\\
\;0           & b^{\,2}\,\\
\end{pmatrix}
\,=\, \dfrac{11}{2}\, \left\|\,M\,\right\|^{\,2}.
\end{align*}
Therefore, \,$(\,\mathcal{F},\, \mathcal{G}\,)$\, is a tight continuous $K$-biframe for \,$H$\, with bound \,$11 \,/\, 2$. 
\end{example}

\begin{theorem}
The pair \,$(\,\mathcal{F},\, \mathcal{G}\,)$\, is a continuous $K$-biframe for \,$H$\, with respect to \,$\left(\,\Omega,\, \mu\,\right)$\, if and only if  \,$(\,\mathcal{G},\, \mathcal{F}\,)$\, is a continuous $K$-biframe for \,$H$\, with respect to \,$\left(\,\Omega,\, \mu\,\right)$.
\end{theorem}

\begin{proof}
Let \,$(\,\mathcal{F},\, \mathcal{G}\,)$\, be a continuous $K$-biframe for \,$H$\, with bounds \,$A$\, and \,$B$.\,Then for each \,$f \,\in\, H$, we have 
\[A\,\left\|\,K^{\,\ast}\,f\,\right\|^{\,2} \leq \int\limits_{\Omega}\,\left<\,f,\, F\,(\,w\,)\,\right>\,\left<\,G\,(\,w\,),\, f\,\right>\,d\mu \leq\, B\left\|\,f\,\right\|^{\,2}.\] 
Now, we can write
\begin{align*}
\int\limits_{\Omega}\,\left<\,f,\, F\,(\,w\,)\,\right>\,\left<\,G\,(\,w\,),\, f\,\right>\,d\mu &\,=\, \overline{\int\limits_{\Omega}\,\left<\,f,\, F\,(\,w\,)\,\right>\,\left<\,G\,(\,w\,),\, f\,\right>\,d\mu}\\
&=\, \int\limits_{\Omega}\,\left<\,f,\, G\,(\,w\,)\,\right>\,\left<\,F\,(\,w\,),\, f\,\right>\,d\mu.
\end{align*}
Thus, for each \,$f \,\in\, H$, we have 
\[A\,\left\|\,K^{\,\ast}\,f\,\right\|^{\,2} \leq \int\limits_{\Omega}\,\left<\,f,\, G\,(\,w\,)\,\right>\,\left<\,F\,(\,w\,),\, f\,\right>\,d\mu \leq\, B\left\|\,f\,\right\|^{\,2}.\]
Therefore, \,$(\,\mathcal{G},\, \mathcal{F}\,)$\, is a continuous $K$-biframe for \,$H$.

Similarly, we can prove the converse part of this Theorem.  
\end{proof}

Next Theorem gives a characterization of continuous $K$-biframe with respect to its frame operator. 

\begin{theorem}\label{3.thm3.8}
Let \,$(\,\mathcal{F},\, \mathcal{G}\,)$\, be a continuous biframe Bessel mapping for \,$H$\, and \,$K \,\in\, \mathcal{B}(\,H\,)$.\,Then \,$(\,\mathcal{F},\, \mathcal{G}\,)$\, is a continuous $K$-biframe for \,$H$\, if and only if there exists a positive constant \,$A$\, such that \,$S_{\mathcal{F},\,\mathcal{G}} \,\geq\, A\,K\,K^{\,\ast}$, where \,$S_{\mathcal{F},\,\mathcal{G}}$\, is the frame operator of \,$(\,\mathcal{F},\, \mathcal{G}\,)$. 
\end{theorem}

\begin{proof}
First, we suppose that \,$(\,\mathcal{F},\, \mathcal{G}\,)$\, be a continuous $K$-biframe for \,$H$\, with lower bound \,$A$\, and  \,$S_{\mathcal{F},\,\mathcal{G}}$\, be its frame operator.\,Then for each \,$f \,\in\, H$, we have 
\[A\,\left\|\,K^{\,\ast}\,f\,\right\|^{\,2} \leq \int\limits_{\Omega}\,\left<\,f,\, F\,(\,w\,)\,\right>\,\left<\,G\,(\,w\,),\, f\,\right>\,d\mu \,=\, \left<\,S_{\mathcal{F},\,\mathcal{G}}\,f,\, f\,\right>.\]
This implies that \,$A\,\left<\,K\,K^{\ast}\,f,\, f\,\right> \,\leq\, \left<\,S_{\mathcal{F},\,\mathcal{G}}\,f,\, f\,\right>$, for all \,$f \,\in\, H$\, and from Theorem \ref{2.5th1.001}, it follows that \,$S_{\mathcal{F},\,\mathcal{G}} \,\geq\, A\,K\,K^{\,\ast}$.

Converse part is obvious.     
\end{proof}

In the following Theorem, we will see that every continuous biframe is a continuous $K$-biframe for \,$H$.

\begin{theorem}\label{3.thm3.9}
Let \,$K \,\in\, \mathcal{B}(\,H\,)$.\,Then the following statements hold:
\begin{itemize}
\item[$(i)$]Every continuous biframe for \,$H$\, is a  continuous $K$-biframe for \,$H$. 
\item[$(ii)$]If \,$\mathcal{R}(\,K\,)$\, is closed, every continuous $K$-biframe for \,$H$\, is a  continuous biframe for \,$\mathcal{R}(\,K\,)$.
\end{itemize} 
\end{theorem}

\begin{proof}$(i)$
Let \,$(\,\mathcal{F},\, \mathcal{G}\,)$\, be a continuous biframe for \,$H$\, with bounds \,$A$\, and \,$B$.\,Then for each \,$f \,\in\, H$, we have 
\begin{align*}
\dfrac{A}{\left\|\,K\,\right\|^{2}}\,\left\|\,K^{\,\ast}\,f\,\right\|^{\,2} \,\leq\, A\,\left\|\,f\,\right\|^{\,2} \leq \int\limits_{\Omega}\,\left<\,f,\, F\,(\,w\,)\,\right>\,\left<\,G\,(\,w\,),\, f\,\right>\,d\mu \leq\, B\left\|\,f\,\right\|^{\,2}.
\end{align*}
Hence, \,$(\,\mathcal{F},\, \mathcal{G}\,)$\, is a continuous $K$-biframe for \,$H$\, with bounds \,$A \,/\, \left\|\,K\,\right\|^{2}$\, and \,$B$.

Proof of $(ii)$.\,\,Let \,$(\,\mathcal{F},\, \mathcal{G}\,)$\, be a continuous $K$-biframe for \,$H$\, with bounds \,$A$\, and \,$B$.\,Since \,$\mathcal{R}(\,K\,)$\, is closed, by Theorem \ref{thm1.1}, there exists a bounded linear operator \,$K^{\dagger}$\, such that \,$K\,K^{\dagger}\,f \,=\, f$,  for all \,$f \,\in\, \mathcal{R}(\,K\,)$.\,Then for each \,$f \,\in\, \mathcal{R}(\,K\,)$, we have 
\begin{align*}
\dfrac{A}{\left\|\,K^{\dagger}\,\right\|^{2}}\,\left\|\,f\,\right\|^{\,2} \,\leq\, A\,\left\|\,K^{\,\ast}\,f\,\right\|^{\,2} \leq \int\limits_{\Omega}\,\left<\,f,\, F\,(\,w\,)\,\right>\,\left<\,G\,(\,w\,),\, f\,\right>\,d\mu \leq\, B\left\|\,f\,\right\|^{\,2}.
\end{align*} 
Thus, \,$(\,\mathcal{F},\, \mathcal{G}\,)$\, is a continuous biframe for \,$\mathcal{R}(\,K\,)$\, with bounds \,$A \,/\, \left\|\,K^{\dagger}\,\right\|^{2}$\, and \,$B$.\,This completes the proof. 
\end{proof}

Now, we validate Theorem \ref{3.thm3.9} by an example. 

\begin{example}
Let \,$H \,=\, \mathbb{R}^{\,3}$.Suppose
\begin{align*}
&\left\{\,f_{\,i}\,\right\}_{i \,=\, 1}^{\,3} \,=\, \left\{\,(\,5,\, 1,\, 1\,),\, (\,-\, 1,\, 7,\, -\, 1\,),\, (\,-\, 1,\, 1,\, 11\,)\,\right\}\,,\\
&\left\{\,g_{\,i}\,\right\}_{i \,=\, 1}^{\,3} \,=\, \left\{\,(\,1,\, 0,\, 0\,),\, (\,0,\, 1,\, 0\,),\, (\,0,\, 0,\, 1\,)\,\right\}.
\end{align*}
Consider \,$\Omega \,=\, \left\{\,x \,\in\, \mathbb{R}^{\,3} \,:\, \|\,x\,\| \,\leq\, 1\,\right\}$.\,Then \,$\left(\,\Omega,\, \mu\,\right)$\, is a measure space, where \,$\mu$\, is the Lebesgue measure.\,Suppose \,$\{\,B_{\,1},\,B_{\,2},\, B_{\,3}\,\}$\, is a partition of \,$\Omega$\, where \,$\mu\,(\,B_{1}\,) \,\geq\, \mu\,(\,B_{2}\,) \,\geq\, \mu\,(\,B_{3}\,) \,>\, 1$.\,Now, we define 
\[F \,:\, \Omega \,\to\, H\hspace{.5cm}\text{by} \hspace{.3cm} F\,(\,w\,) \,=\, \begin{cases}
\dfrac{f_{1}}{\sqrt{\,\mu\left(\,B_{1}\,\right)}} & \text{if\;\;}\; w \,\in\, B_{1} \\ \dfrac{f_{2}}{\sqrt{\,\mu\left(\,B_{2}\,\right)}} & \text{if\;\;}\; w \,\in\, B_{2}\\ \dfrac{f_{3}}{\sqrt{\,\mu\left(\,B_{3}\,\right)}} & \text{if\;\;}\; w \,\in\, B_{3}\,, \end{cases}\]
\[G \,:\, \Omega \,\to\, H\hspace{.5cm}\text{by} \hspace{.3cm} G\,(\,w\,) \,=\, \begin{cases}
\dfrac{g_{1}}{\sqrt{\,\mu\left(\,B_{1}\,\right)}} & \text{if\;\;}\; w \,\in\, B_{1} \\ \dfrac{g_{2}}{\sqrt{\,\mu\left(\,B_{2}\,\right)}} & \text{if\;\;}\; w \,\in\, B_{2}\\ \dfrac{g_{3}}{\sqrt{\,\mu\left(\,B_{3}\,\right)}} & \text{if\;\;}\; w \,\in\, B_{3} \end{cases}\]
It is easy to verify that for all \,$f \,\in\, H$, \,$w \,\mapsto\, \left<\,f,\, F\,(\,w\,)\,\right>$\, and \,$w \,\mapsto\, \left<\,f,\, G\,(\,w\,)\,\right>$\, are measurable functions on \,$\Omega$.\,Now, for \,$f \,=\, \left(\,x,\, y,\, z\,\right) \,\in\, H$, we have
\begin{align*}
&\int\limits_{\,\Omega}\,\left<\,f,\, F\,(\,w\,)\,\right>\,\left<\,G\,(\,w\,),\, f\,\right>\,d\mu \\
&= \int\limits_{\,B_{1}}\left<\,f,\, \dfrac{f_{\,1}}{\sqrt{\mu\left(\,B_{1}\,\right)}}\,\right>\,\left<\,\dfrac{g_{\,1}}{\sqrt{\mu\left(\,B_{1}\,\right)}},\, f\,\right>\,d\mu\\
& \,+\, \int\limits_{\,B_{2}}\left<\,f,\, \dfrac{f_{\,2}\,}{\sqrt{\mu\left(\,B_{2}\,\right)}}\,\right>\,\left<\,\dfrac{\, g_{\,2}}{\sqrt{\mu\left(\,B_{2}\,\right)}},\, f\,\right>\,d\mu\\
&+\,\int\limits_{\,B_{3}}\left<\,f,\, \dfrac{f_{\,3}\,}{\sqrt{\mu\left(\,B_{3}\,\right)}}\,\right>\,\left<\, \dfrac{g_{\,3}\,}{\sqrt{\mu\left(\,B_{3}\,\right)}},\, f\,\right>\,d\mu\\
&=\,\left<\,f,\, f_{1}\,\right>\,\left<\,g_{1},\, f\,\right> \,+\, \left<\,f,\, f_{2}\,\right>\,\left<\,g_{2},\, f\,\right> \,+\, \left<\,f,\, f_{3}\,\right>\,\left<\,g_{3},\, f\,\right>\\
&=\,\left<\,\left(\,x,\, y,\, z\,\right),\, (\,5,\, 1,\, 1\,)\,\right>\,\left<\,(\,1,\, 0,\, 0\,),\, \left(\,x,\, y,\, z\,\right)\,\right> \,+\,\\
&\hspace{1cm}+\, \left<\,\left(\,x,\, y,\, z\,\right),\, (\,-\, 1,\, 7,\, \,-\, 1\,)\,\right>\,\left<\,(\,0,\, 1,\, 0\,),\, \left(\,x,\, y,\, z\,\right)\,\right> \\
&\hspace{1cm}+\,\left<\,\left(\,x,\, y,\, z\,\right),\, (\,-\, 1,\, 1,\, 11\,)\,\right>\,\left<\,(\,0,\, 0,\, 1\,),\, \left(\,x,\, y,\, z\,\right)\,\right>\\
&=\,\left(\,5\,x \,+\, y \,+\, z\,\right)\,x \,+\, \left(\,-\, x \,+\, 7\,y \,-\, z\,\right)\,y \,+\, \left(\,-\, x \,+\, y \,+\, 11\,z\,\right)\,z\\
&=\,5\,x^{\,2} \,+\, 7\,y^{\,2} \,+\, 11\,z^{\,2} \,\leq\, 11\,\left(\,x^{\,2} \,+\, y^{\,2} \,+\, z^{\,2}\,\right) \,=\, 11\,\left\|\,\left(\,x,\, y,\, z\,\right)\,\right\|^{2} \,=\, 11\,\left\|\,f\,\right\|^{2}. 
\end{align*} 
Therefore, \,$(\,F,\, G\,)$\, is a continuous biframe for \,$H$\, with bounds \,$5$\, and \,$11$.\,Define \,$K \,:\, \mathbb{R}^{\,3} \,\to\, \mathbb{R}^{\,3}$\, by \,$K\left(\,x,\, y,\, z\,\right) \,=\, \left(\,2\,x,\,  \,-\, 2\,y,\, \,-\, 2\,z\,\right)$, for every \,$\left(\,x,\, y,\, z\,\right) \,\in\, \mathbb{R}^{\,3}$.\,Then \,$K$\, is bounded linear operator on \,$\mathbb{R}^{\,3}$.\,The matrix associated with the operator \,$K$\, is given by
\[
\left[\,K\,\right] \,=\, 
\begin{pmatrix}
\;2 & 0 & 0\\
\; 0   & \,-\, 2 & \,0\\
\; 0    & 0  & \,-\, 2
\end{pmatrix}
.\] 
Thus, for every \,$\left(\,x,\, y,\, z\,\right) \,\in\, \mathbb{R}^{\,3}$, we get \,$K^{\,\ast}\left(\,x,\, y,\, z\,\right) \,=\, \left(\,2\,x,\,  \,-\, 2\,y,\, \,-\, 2\,z\,\right)$\, and \,$\left\|\,K^{\,\ast}\,\left(\,x,\, y,\, z\,\right)\,\right\|^{\,2} \,=\, 4\,\left(\,x^{\,2} \,+\, y^{\,2} \,+\, z^{\,2}\,\right)$.\,Thus, for each \,$f \,\in\, H$, we get
\begin{align*}
\dfrac{5}{4}\,\left\|\,K^{\,\ast}\,f\,\right\|^{2} \,\leq\, \int\limits_{\,\Omega}\,\left<\,f,\, F\,(\,w\,)\,\right>\,\left<\,G\,(\,w\,),\, f\,\right>\,d\mu \,\leq\, 11\,\left\|\,f\,\right\|^{2}.
\end{align*}
Hence, \,$(\,F,\, G\,)$\, is a continuous $K$-biframe for \,$H$\, with bounds \,$5 \,/\, 4$\, and \,$11$. 
\end{example}

\begin{theorem}\label{3.tmn2.19}
Let \,$(\,\mathcal{F},\, \mathcal{G}\,)$\, be a continuous $K_{1}$-biframe and continuous $K_{2}$-biframe for \,$H$.\,Then for any scalar \,$\alpha,\, \beta$, \,$(\,\mathcal{F},\, \mathcal{G}\,)$\, is a continuous $\alpha\,K_{1} \,+\, \beta\,K_{2}$-biframe and continuous \,$K_{1}\,K_{2}$-biframe for \,$H$.  
\end{theorem}

\begin{proof}
Since \,$(\,\mathcal{F},\, \mathcal{G}\,)$\, is a continuous \,$K_{1}$-biframe and continuous \,$K_{2}$-biframe for \,$H$, there exist positive constants \,$A_{n},\, B_{n} \,>\, 0$, \,$n \,=\, 1,\, 2$\, such that
\begin{align}
A_{n}\,\left \|\,K_{n}^{\,\ast}\,f\,\right \|^{\,2}& \,\leq\, \int\limits_{\Omega}\,\left<\,f,\, F\,(\,w\,)\,\right>\,\left<\,G\,(\,w\,),\, f\,\right>\,d\mu \,\leq\, B_{n}\,\left\|\,f\,\right\|^{\,2},\label{3.emn2.21}
\end{align}
for all \,$f \,\in\, H$.\,Now, for each \,$f \,\in\, H$, we have
\begin{align*}
&\left \|\,K_{1}^{\,\ast}\,f\,\right \|^{\,2} \,=\, \dfrac{1}{|\,\alpha\,|^{\,2}}\,\left \|\,\alpha\,K_{1}^{\,\ast}\,f\,\right \|^{\,2}\\
&=\,\dfrac{1}{|\,\alpha\,|^{\,2}}\,\left \|\,\left(\,\alpha\,K_{1}^{\,\ast} \,+\, \beta\,K_{2}^{\,\ast}\,\right)\,f \,-\, \beta\,K_{2}^{\,\ast}\,f\,\right \|^{\,2}\\
&\geq\,\dfrac{1}{|\,\alpha\,|^{\,2}}\,\left(\,\left \|\,\left(\,\alpha\,K_{1}^{\,\ast} \,+\, \beta\,K_{2}^{\,\ast}\,\right)\,f\,\right \|^{\,2} \,-\, \left \|\,\beta\,K_{2}^{\,\ast}\,f\,\right \|^{\,2}\,\right).\\
&\Rightarrow\,\left \|\,\left(\,\alpha\,K_{1}^{\,\ast} \,+\, \beta\,K_{2}^{\,\ast}\,\right)\,f\,\right \|^{\,2} \\
&\leq\,|\,\alpha\,|^{\,2}\,\left \|\,K_{1}^{\,\ast}\,f\,\right \|^{\,2} \,+\, |\,\beta\,|^{\,2}\,\left \|\,K_{2}^{\,\ast}\,f\,\right \|^{\,2}\\
&\leq\,\max\left\{\,|\,\alpha\,|^{\,2},\, |\,\beta\,|^{\,2}\,\right\}\left\{\left \|\,K_{1}^{\,\ast}\,f\,\right \|^{\,2} \,+\, \left \|\,K_{2}^{\,\ast}\,f\,\right \|^{\,2}\right\}. 
\end{align*}
Thus, for each \,$f \,\in\, H$, we have
\begin{align*}
&\dfrac{1}{\max\left\{\,|\,\alpha\,|^{\,2},\, |\,\beta\,|^{\,2}\,\right\}}\,\left \|\,\left(\,\alpha\,K_{1}^{\,\ast} \,+\, \beta\,K_{2}^{\,\ast}\,\right)\,f\,\right \|^{\,2} \\
&\leq\, \left(\,\dfrac{1}{A_{1}} \,+\, \dfrac{1}{A_{2}}\,\right)\,\int\limits_{\Omega}\,\left<\,f,\, F\,(\,w\,)\,\right>\,\left<\,G\,(\,w\,),\, f\,\right>\,d\mu.
\end{align*} 
On the other hand, from (\ref{3.emn2.21}), for each \,$f \,\in\, H$, we get
\[\int\limits_{\Omega}\,\left<\,f,\, F\,(\,w\,)\,\right>\,\left<\,G\,(\,w\,),\, f\,\right>\,d\mu \,\leq\, \dfrac{B_{1} \,+\, B_{2}}{2}\,\left \|\,f\,\right \|^{\,2}.\]
Hence, \,$(\,\mathcal{F},\, \mathcal{G}\,)$\, is a continuous $\alpha\,K_{1} \,+\, \beta\,K_{2}$-biframe for \,$H$.\,Moreover, for each \,$f \,\in\, H$, we have
\begin{align*}
\left \|\,\left(\,K_{1}\,K_{2}\,\right)^{\,\ast}\,f\,\right \|^{\,2} &\,=\, \left \|\,K_{2}^{\,\ast}\,K_{1}^{\,\ast}\,f\,\right \|^{\,2}\leq\,\left\|\,K_{2}^{\,\ast}\,\right\|^{\,2}\,\left \|\,K_{1}^{\,\ast}\,f\,\right \|^{\,2}.
\end{align*} 
Thus, for each \,$f \,\in\, H$, we have
\begin{align*}
\dfrac{A_{1}}{\left\|\,K_{2}\,\right\|^{\,2}}\,\left \|\,\left(\,K_{1}\,K_{2}\,\right)^{\,\ast}\,f\,\right \|^{\,2} \,\leq\, A_{1}\,\left \|\,K_{1}^{\,\ast}\,f\,\right \|^{\,2} & \leq\,\int\limits_{\Omega}\,\left<\,f,\, F\,(\,w\,)\,\right>\,\left<\,G\,(\,w\,),\, f\,\right>\,d\mu\,\\
& \leq B_{1}\,\left \|\,f\,\right \|^{\,2}.
\end{align*}
Therefore, \,$(\,\mathcal{F},\, \mathcal{G}\,)$\, is a continuous \,$K_{1}\,K_{2}$-biframe for \,$H$.\,This completes the proof.      
\end{proof}

Next Theorem is a generalization of the Theorem \ref{3.tmn2.19}. 

\begin{theorem}
Let \,$K_{j} \,\in\, \mathcal{B}\,\left(\,H\,\right), \,\, j \,=\, 1,\, 2,\, \cdots,\, n$\, and \,$(\,\mathcal{F},\, \mathcal{G}\,)$\, be a continuous \,$K_{j}$-biframe for \,$H$.\,Then
\begin{itemize}
\item[$(i)$]If \,$a_{\,j}$, for \,$j \,=\, 1,\, 2,\, \cdots,\, n$,  are finite collection of scalars, then \,$(\,\mathcal{F},\, \mathcal{G}\,)$\, is a continuous \,$\sum\limits^{n}_{j \,=\, 1}\,a_{\,j}\,K_{j}$-biframe for \,$H$.
\item[$(ii)$]\,\,\,$(\,\mathcal{F},\, \mathcal{G}\,)$\, is a continuous \,$\displaystyle \prod_{j \,=\, 1}^{\,n}\,K_{j}$-biframe for \,$H$. 
\end{itemize}
\end{theorem}

\begin{proof}$(i)$\,\,
Since \,$(\,\mathcal{F},\, \mathcal{G}\,)$\, is a continuous \,$K_{j}$-biframe for \,$H$, for all \,$j$, there exist constants \,$A,\,B \,>\, 0$\, such that
\begin{align*}
A\,\left \|\,K_{j}^{\,\ast}\,f\,\right \|^{\,2}& \,\leq\, \int\limits_{\Omega}\,\left<\,f,\, F\,(\,w\,)\,\right>\,\left<\,G\,(\,w\,),\, f\,\right>\,d\mu \,\leq\, B\,\left\|\,f\,\right\|^{\,2}.
\end{align*}
Then for each \,$f \,\in\, H$, we have
\begin{align*}
&\dfrac{A}{n\,\max\limits_{j}\,|\,a_{\,j}\,|^{\,2}}\,\left\|\,\left(\,\sum\limits^{n}_{j \,=\, 1}\,a_{\,j}\,K_{\,j}\,\right)^{\,\ast}\,f\,\right\|^{\,2} \,\leq\, A\,\left\|\,K_{j}^{\,\ast}\,f\,\right\|^{\,2}\\
&\,\leq\, \int\limits_{\Omega}\,\left<\,f,\, F\,(\,w\,)\,\right>\,\left<\,G\,(\,w\,),\, f\,\right>\,d\mu \,\leq\, B\,\left\|\,f\,\right\|^{\,2}.
\end{align*}
Thus, \,$(\,\mathcal{F},\, \mathcal{G}\,)$\, is a continuous \,$\sum\limits^{n}_{j \,=\, 1}\,a_{\,j}\,K_{j}$-biframe for \,$H$.

Proof of \,$(ii)$.\,\, For each \,$f \,\in\, H$, we have
\begin{align*}
&\dfrac{A}{\displaystyle \prod_{j \,=\, 1}^{n \,-\, 1}\,\left\|\,K^{\,\ast}_{j}\,\right\|^{\,2}}\,\left\|\,\left(\,\displaystyle \prod_{j \,=\, 1}^{\,n}\,K_{j}\,\right)^{\,\ast}\,f\,\right\|^{\,2} \,\leq\, A\,\left\|\,K^{\,\ast}_{n}\,f\,\right\|^{\,2}\\
&\,\leq\, \int\limits_{\Omega}\,\left<\,f,\, F\,(\,w\,)\,\right>\,\left<\,G\,(\,w\,),\, f\,\right>\,d\mu \,\leq\, B\,\left\|\,f\,\right\|^{\,2}.
\end{align*} 
Thus, \,$(\,\mathcal{F},\, \mathcal{G}\,)$\, is a continuous \,$\displaystyle \prod_{j \,=\, 1}^{\,n}\,K_{j}$-biframe for \,$H$.\,This completes the proof.  
\end{proof}

\begin{theorem}
Let \,$K \,\in\, \mathcal{B}(\,H\,)$\, and \,$A_{1},\, A_{2} \,>\, 0$.\,Suppose \,$(\,\mathcal{F},\, \mathcal{G}\,)$\, is a tight continuous $K$-biframe for \,$H$\, with bound \,$A_{1}$.\,Then \,$(\,\mathcal{F},\, \mathcal{G}\,)$\, is a tight continuous biframe for \,$H$\, with bound \,$A_{2}$\, if and only if \,$K\,\left(\,A_{1} \,/\, A_{2}\,\right)\,K^{\,\ast} \,=\, I_{H}$.   
\end{theorem}

\begin{proof}
First, we suppose that \,$(\,\mathcal{F},\, \mathcal{G}\,)$\, is a tight continuous biframe for \,$H$\, with bound \,$A_{2}$.\,Then for each \,$f \,\in\, H$, we have
\begin{align*}
\int\limits_{\Omega}\,\left<\,f,\, F\,(\,w\,)\,\right>\,\left<\,G\,(\,w\,),\, f\,\right>\,d\mu \,=\, A_{2}\,\left \|\,f\,\right \|^{\,2}.
\end{align*}
Since \,$(\,\mathcal{F},\, \mathcal{G}\,)$\, is a tight continuous $K$-biframe for \,$H$\, with bound \,$A_{1}$, for each \,$f \,\in\, H$, we have  
\begin{align*}
\int\limits_{\Omega}\,\left<\,f,\, F\,(\,w\,)\,\right>\,\left<\,G\,(\,w\,),\, f\,\right>\,d\mu \,=\, A_{1}\,\left \|\,K^{\,\ast}\,f\,\right \|^{\,2}.
\end{align*}
From the above two equalities, for each \,$f \,\in\, H$, we get
\begin{align*}
\left \|\,K^{\,\ast}\,f\,\right \|^{\,2} \,=\, \left(\,A_{2} \,/\, A_{1}\,\right)\,\left \|\,f\,\right \|^{\,2} \,\Rightarrow\, \left<\,K\,K^{\,\ast}\,f,\, f\,\right> \,=\, \left<\,\left(\,A_{2} \,/\, A_{1}\,\right)\,f,\, f\,\right>.
\end{align*}
This implies that \,$K\,\left(\,A_{1} \,/\, A_{2}\,\right)\,K^{\,\ast} \,=\, I_{H}$. 

Next, we suppose that \,$K\,\left(\,A_{2} \,/\, A_{1}\,\right)\,K^{\,\ast} \,=\, I_{H}$.\,This gives that \,$\left \|\,K^{\,\ast}\,f\,\right \|^{\,2} \,=\, \left(\,A_{2} \,/\, A_{1}\,\right)\,\left \|\,f\,\right \|^{\,2}$.\,Now, \,$(\,\mathcal{F},\, \mathcal{G}\,)$\, is a tight continuous $K$-biframe for \,$H$\, with bound \,$A_{1}$, so \,$(\,\mathcal{F},\, \mathcal{G}\,)$\, is a tight continuous biframe for \,$H$\, with bound \,$A_{2}$.\,This completes the proof.    
\end{proof}

\begin{corollary}
Let \,$K \,\in\, \mathcal{B}(\,H\,)$.\,Suppose \,$(\,\mathcal{F},\, \mathcal{G}\,)$\, is a Parseval continuous $K$-biframe for \,$H$.\,Then \,$(\,\mathcal{F},\, \mathcal{G}\,)$\, is a Parseval continuous $K$-biframe for \,$H$\, if and only if \,$K\,K^{\,\ast} \,=\, I_{H}$.
\end{corollary}

\section{Continuous $K$-biframes and operators}

In this section, we construct continuous $K$-biframe by using a bounded linear operator on \,$H$\, and finally establish a relationship between quotient operator and continuous $K$-biframe for \,$H$.

\begin{theorem}
Let \,$K \,\in\, \mathcal{B}(\,H\,)$\, and \,$(\,\mathcal{F},\, \mathcal{G}\,)$\, be a continuous biframe for \,$H$.\,Then \,$(\,K\,\mathcal{F},\, K\,\mathcal{G}\,)$\, is a continuous $K$-biframe for \,$H$.
\end{theorem}

\begin{proof}
Suppose \,$(\,\mathcal{F},\, \mathcal{G}\,)$\, is a continuous biframe for \,$H$\, with bounds \,$A$\, and \,$B$.\,Now, for each \,$f \,\in\, H$, we have
\begin{align*}
\int\limits_{\Omega}\,\left<\,f,\, K\,F\,(\,w\,)\,\right>\,\left<\,K\,G\,(\,w\,),\, f\,\right>\,d\mu &\,=\, \int\limits_{\Omega}\,\left<\,K^{\,\ast}\,f,\, F\,(\,w\,)\,\right>\,\left<\,G\,(\,w\,),\, K^{\,\ast}\,f\,\right>\,d\mu\\
&\geq\,A\,\left \|\,K^{\,\ast}\,f\,\right \|^{\,2}.
\end{align*}
On the other hand 
\begin{align*}
\int\limits_{\Omega}\,\left<\,f,\, K\,F\,(\,w\,)\,\right>\,\left<\,K\,G\,(\,w\,),\, f\,\right>\,d\mu &\,\leq\, B\,\left \|\,K^{\,\ast}\,f\,\right \|^{\,2}\\
&\leq\,B\,\left\|\,K\,\right\|^{\,2}\,\left \|\,f\,\right \|^{\,2}.
\end{align*} 
This completes the proof.
\end{proof}

The frame operator \,$S^{\,\prime}_{F,\,G}$\, is given by
\begin{align*}
S^{\,\prime}_{F,\,G}\,f &\,=\, \int\limits_{\,\Omega}\,\left<\,f, K\,F\,(\,w\,)\,\right>\,K\,G\,(\,w\,)\,d\mu\\
&=\,K\,\left(\,\int\limits_{\,\Omega}\,\left<\,K^{\,\ast}\,f, F\,(\,w\,)\,\right>\,G\,(\,w\,)\,d\mu\,\right) \,=\, K\,S_{F,\,G}\,K^{\,\ast}\,f,
\end{align*}
for all \,$f \,\in\, H$.

\begin{theorem}
Let \,$K \,\in\, \mathcal{B}(\,H\,)$.\,Suppose \,$(\,\mathcal{F},\, \mathcal{G}\,)$\, is a continuous biframe for \,$H$\, with frame operator \,$S_{\mathcal{F},\,\mathcal{G}}$.\,Then \,$\left(\,K\,S^{\,-\, 1}_{\mathcal{F},\,\mathcal{G}}\,\mathcal{F},\, K\,S^{\,-\, 1}_{\mathcal{F},\,\mathcal{G}}\,\mathcal{G}\,\right)$\, is a continuous $K$-biframe for \,$H$.
\end{theorem}

\begin{proof}
Let \,$(\,\mathcal{F},\, \mathcal{G}\,)$\, be a continuous biframe for \,$H$\, with bounds \,$A,\,B$\, and \,$S_{\mathcal{F},\,\mathcal{G}}$\, be the corresponding frame operator.\,Since \,$S_{\mathcal{F},\,\mathcal{G}}$\, is invertible, for \,$f \,\in\, H$, we get
\begin{align*}
\left\|\,f\,\right\|^{\,2} \,\leq\, \left\|\,S_{\mathcal{F},\,\mathcal{G}}\,\right\|^{\,2}\,\left\|\,S^{\,-\, 1}_{\mathcal{F},\,\mathcal{G}}\,f\,\right\|^{\,2} \,\Rightarrow\, \left\|\,K^{\,\ast}\,f\,\right\|^{\,2} \,\leq\, \left\|\,S_{\mathcal{F},\,\mathcal{G}}\,\right\|^{\,2}\,\left\|\,S^{\,-\, 1}_{\mathcal{F},\,\mathcal{G}}\,K^{\,\ast}\,f\,\right\|^{\,2}.   
\end{align*}
Now, for each \,$f \,\in\, H$, we have
\begin{align*}
&\int\limits_{\Omega}\,\left<\,f,\, K\,S^{\,-\, 1}_{\mathcal{F},\,\mathcal{G}}\,F\,(\,w\,)\,\right>\,\left<\,K\,S^{\,-\, 1}_{\mathcal{F},\,\mathcal{G}}\,G\,(\,w\,),\, f\,\right>\,d\mu \\
&\,=\, \int\limits_{\Omega}\,\left<\,S^{\,-\, 1}_{\mathcal{F},\,\mathcal{G}}\,K^{\,\ast}\,f,\, F\,(\,w\,)\,\right>\,\left<\,G\,(\,w\,),\, S^{\,-\, 1}_{\mathcal{F},\,\mathcal{G}}\,K^{\,\ast}\,f\,\right>\,d\mu\\
&\geq\,A\,\left \|\,S^{\,-\, 1}_{\mathcal{F},\,\mathcal{G}}\,K^{\,\ast}\,f\,\right \|^{\,2} \,\geq\, \dfrac{A}{\left\|\,S_{\mathcal{F},\,\mathcal{G}}\,\right\|^{\,2}}\,\left \|\,K^{\,\ast}\,f\,\right \|^{\,2}.
\end{align*}
On the other hand
\begin{align*}
&\int\limits_{\Omega}\,\left<\,f,\, K\,S^{\,-\, 1}_{\mathcal{F},\,\mathcal{G}}\,F\,(\,w\,)\,\right>\,\left<\,K\,S^{\,-\, 1}_{\mathcal{F},\,\mathcal{G}}\,G\,(\,w\,),\, f\,\right>\,d\mu \\
&\leq\,B\,\left \|\,S^{\,-\, 1}_{\mathcal{F},\,\mathcal{G}}\,K^{\,\ast}\,f\,\right \|^{\,2} \,\leq\,B\,\left\|\,S^{\,-\, 1}_{\mathcal{F},\,\mathcal{G}}\,\right\|^{\,2}\,\left \|\,K\,\right \|^{\,2}\,\left \|\,f\,\right \|^{\,2}.
\end{align*}
Thus, \,$\left(\,K\,S^{\,-\, 1}_{\mathcal{F},\,\mathcal{G}}\,\mathcal{F},\, K\,S^{\,-\, 1}_{\mathcal{F},\,\mathcal{G}}\,\mathcal{G}\,\right)$\, is a continuous $K$-biframe for \,$H$.\,This completes the proof.
\end{proof}

\begin{theorem}\label{3.ttm3.20}
Let \,$(\,\mathcal{F},\, \mathcal{G}\,)$\, be a continuous $K$-biframe for \,$H$\, and \,$U \,\in\, \mathcal{B}\left(\,H\,\right)$.\,Then \,$\left(\,U\,\mathcal{F},\, U\,\mathcal{G}\,\right)$\, is a continuous $U\,K$-biframe for \,$H$.
\end{theorem}

\begin{proof}
Now, for each \,$f \,\in\, H$, we have
\begin{align*}
\int\limits_{\Omega}\,\left<\,f,\, U\,F\,(\,w\,)\,\right>\,\left<\,U\,G\,(\,w\,),\, f\,\right>\,d\mu &\,=\, \int\limits_{\Omega}\,\left<\,U^{\,\ast}\,f,\, F\,(\,w\,)\,\right>\,\left<\,G\,(\,w\,),\, U^{\,\ast}\,f\,\right>\,d\mu\\
&\geq\,A\,\left \|\,K^{\,\ast}\,U^{\,\ast}\,f\,\right \|^{\,2}.
\end{align*}
On the other hand
\begin{align*}
\int\limits_{\Omega}\,\left<\,f,\, U\,F\,(\,w\,)\,\right>\,\left<\,U\,G\,(\,w\,),\, f\,\right>\,d\mu \,\leq\, B\,\|\,U\,\|^{\,2}\,\left\|\,f\,\right \|^{\,2}.
\end{align*}
Hence, \,$\left(\,U\,\mathcal{F},\, U\,\mathcal{G}\,\right)$\, is a continuous $U\,K$-biframe for \,$H$. 
\end{proof}

In the next two Theorems, we will construct new type of continuous $K$-biframe from a given continuous $K$-biframe for \,$H$\, by using a bounded linear operator.

\begin{theorem}
Let \,$U \,\in\, \mathcal{B}\,\left(\,H\,\right)$\, be an operator and \,$(\,\mathcal{F},\, \mathcal{G}\,)$\, be a continuous $K$-biframe for \,$H$.\,Then \,$\left(\,U\,\mathcal{F},\, U\,\mathcal{G}\,\right)$\, is a continuous $U\,K\,U^{\,\ast}$-biframe for \,$H$.
\end{theorem}

\begin{proof}
Let \,$(\,\mathcal{F},\, \mathcal{G}\,)$\, be a continuous $K$-biframe for \,$H$\, with bounds \,$A$\, and \,$B$.
Now, for each \,$f \,\in\, H$, we have
\begin{align*}
&\dfrac{A}{\|\,U\,\|^{\,2}}\, \left\|\,\left(\,U\,K\,U^{\,\ast}\,\right)^{\,\ast}\,f\,\right\|^{\,2} \,=\, \dfrac{A}{\|\,U\,\|^{\,2}}\, \left\|\,U\,K^{\,\ast}\,U^{\,\ast}\,f\,\right\|^{\,2}\\
& \,\leq\, A\, \left\|\,K^{\,\ast}\,U^{\,\ast}\,f\,\right\|^{\,2}\leq\, \int\limits_{\Omega}\,\left<\,U^{\,\ast}\,f,\, F\,(\,w\,)\,\right>\,\left<\,G\,(\,w\,),\, U^{\,\ast}\,f\,\right>\,d\mu\\
&=\,\int\limits_{\Omega}\,\left<\,f,\, U\,F\,(\,w\,)\,\right>\,\left<\,U\,G\,(\,w\,),\, f\,\right>\,d\mu \leq\, B\,\|\,U\,\|^{\,2}\,\left\|\,f\,\right \|^{\,2}.
\end{align*}
Thus, \,$\left(\,U\,\mathcal{F},\, U\,\mathcal{G}\,\right)$\, is a continuous $U\,K\,U^{\,\ast}$-biframe for \,$H$.
\end{proof}

\begin{theorem}
Let \,$U \,\in\, \mathcal{B}\,\left(\,H\,\right)$\, be an invertible operator and \,$(\,U\,\mathcal{F},\, U\,\mathcal{G}\,)$\, be a continuous $K$-biframe for \,$H$.\,Then \,$\left(\,\mathcal{F},\, \mathcal{G}\,\right)$\, is a continuous $U^{\,-\, 1}\,K\,U$-biframe for \,$H$.
\end{theorem}

\begin{proof}
Let \,$(\,U\,\mathcal{F},\, U\,\mathcal{G}\,)$\, be a continuous $K$-biframe for \,$H$\, with bounds \,$A$\, and \,$B$.\,Now, for each \,$f \,\in\, H$, we have
\begin{align*}
&\dfrac{A}{\|\,U\,\|^{\,2}}\,\left \|\,\left(\,U^{\,-\, 1}\,K\,U\,\right)^{\,\ast}\,f\,\right \|^{\,2} \,=\, \dfrac{A}{\|\,U\,\|^{\,2}}\,\left\|\,U^{\,\ast}\,K^{\,\ast}\,\left(\,U^{\,-\, 1}\,\right)^{\,\ast}\,f\,\right\|^{\,2}\\
& \,\leq\, A\;\left\|\,K^{\,\ast}\,\left(\,U^{\,-\, 1}\,\right)^{\,\ast}\,f\,\right\|^{\,2}\\
&\leq\, \int\limits_{\Omega}\,\left<\,\left(\,U^{\,-\, 1}\,\right)^{\,\ast}\,f,\, U\,F\,(\,w\,)\,\right>\,\left<\,U\,G\,(\,w\,),\, \left(\,U^{\,-\, 1}\,\right)^{\,\ast}\,f\,\right>\,d\mu \\
&=\, \int\limits_{\Omega}\,\left<\,f,\, F\,(\,w\,)\,\right>\,\left<\,G\,(\,w\,),\, f\,\right>\,d\mu.
\end{align*}
On the other hand, for each \,$f \,\in\, H$, we have
\begin{align*}
&\int\limits_{\Omega}\,\left<\,f,\, F\,(\,w\,)\,\right>\,\left<\,G\,(\,w\,),\, f\,\right>\,d\mu  \,=\, \int\limits_{\Omega}\,\left<\,f,\, U^{\,-\, 1}\,U\,F\,(\,w\,)\,\right>\,\left<\,U^{\,-\, 1}\,U\,G\,(\,w\,),\, f\,\right>\,d\mu \\
&\hspace{1cm}=\, \int\limits_{\Omega}\,\left<\,\left(\,U^{\,-\, 1}\,\right)^{\,\ast}\,f,\, U\,F\,(\,w\,)\,\right>\,\left<\,U\,G\,(\,w\,),\, \left(\,U^{\,-\, 1}\,\right)^{\,\ast}\,f\,\right>\,d\mu \\
&\hspace{1cm}\leq\, B \, \left\|\,\left(\,U^{\,-\, 1}\,\right)^{\,\ast}\,f\, \right\|^{\,2}\,\leq\, B\; \left\|\,U^{\,-\, 1}\,\right\|^{\,2}\,\left\|\,f\,\right\|^{\,2}.
\end{align*}
Thus, \,$\left(\,\mathcal{F},\, \mathcal{G}\,\right)$\, is a continuous $U^{\,-\, 1}\,K\,U$-biframe for \,$H$.
\end{proof}

In the following theorem, we give a characterization of continuous $K$-biframe for \,$H$\, with the help of a bounded linear operator.

\begin{theorem}
Let \,$\left(\,\mathcal{F},\, \mathcal{G}\,\right)$\, is a continuous $K$-biframe for \,$H$\, with bounds \,$A,\, B$\, and \,$U \,:\, H \,\to\, H$\; be a bounded linear operator such that \,$\mathcal{R}\,(\,U\,) \,\subset\, \mathcal{R}\,(\,K\,)$.\;Then \,$(\,U\,\mathcal{F},\, U\,\mathcal{G}\,)$\, is a continuous $K$-biframe for \,$H$\, if and only if there exists a \,$\delta \,>\, 0$\, such that \,$\left\|\,U^{\,\ast}\,f\,\right \| \,\geq\, \delta\,\left \|\,K^{\,\ast}\,f\,\right \|$, for all \,$f \,\in\, H$.
\end{theorem}

\begin{proof}
Let  \,$(\,U\,\mathcal{F},\, U\,\mathcal{G}\,)$\, be a continuous $K$-biframe for \,$H$\, with bounds \,$A,\, B$.\,Then for each \,$f \,\in\, H$, we have
\begin{align*}
&A\,\left\|\,K^{\,\ast}\,f\,\right \|^{\,2} \,\leq\, \int\limits_{\Omega}\,\left<\,f,\, U\,F\,(\,w\,)\,\right>\,\left<\,U\,G\,(\,w\,),\, f\,\right>\,d\mu \\
&\,=\, \int\limits_{\Omega}\,\left<\,U^{\,\ast}\,f,\, F\,(\,w\,)\,\right>\,\left<\,G\,(\,w\,),\, U^{\,\ast}\,f\,\right>\,d\mu  \,\leq\, B\,\left\|\,U^{\,\ast}\,f\,\right \|^{\,2}. 
\end{align*}
This implies that for each \,$f \,\in\, H$, we get
\[\left\|\,U^{\,\ast}\,f\,\right \| \,\geq\, \sqrt{A \,/\, B}\,\left\|\,K^{\,\ast}\,f\,\right \|.\] 
Conversely, for each \,$f \,\in\, H$, we can write 
\begin{align*}
\left\|\,U^{\,\ast}\,f\,\right\|^{\,2} &\,=\, \left\|\,(\,K^{\,\dagger}\,)^{\,\ast}\,K^{\,\ast}\,U^{\,\ast}\,f\,\right\|^{\,2} \,\leq\, \left\|\,K^{\,\dagger}\,\right\|^{\,2}\,\left\|\,K^{\,\ast}\,U^{\,\ast}\,f\,\right\|^{\,2}.
\end{align*}
Now, for each \,$f \,\in\, H$, we have
\begin{align*}
&A\,\delta^{\,2}\,\left\|\,K^{\,\dagger}\,\right\|^{\,-\, 2}\,\left\|\,K^{\,\ast}\,f\,\right\|^{\,2} \,\leq\, A\,\left\|\,K^{\,\dagger}\,\right\|^{\,-\, 2}\,\left\|\,U^{\,\ast}\,f\,\right\|^{\,2}\\
&\leq\, A\,\left\|\,K^{\,\ast}\,U^{\,\ast}\,f\,\right\|^{\,2}\leq\, \int\limits_{\Omega}\,\left<\,U^{\,\ast}\,f,\, F\,(\,w\,)\,\right>\,\left<\,G\,(\,w\,),\, U^{\,\ast}\,f\,\right>\,d\mu  \\
&=\,\int\limits_{\Omega}\,\left<\,f,\, U\,F\,(\,w\,)\,\right>\,\left<\,U\,G\,(\,w\,),\, f\,\right>\,d\mu \leq\,B\,\left\|\,U\,\right \|^{\,2}\, \left\|\,f\,\right\|^{\,2}. 
\end{align*}
Thus, \,$(\,U\,\mathcal{F},\, U\,\mathcal{G}\,)$\, is a continuous $K$-biframe for \,$H$.\,This completes the proof. 
\end{proof}

\begin{theorem}\label{3.th3.22}
Let \,$\left(\,\mathcal{F},\, \mathcal{G}\,\right)$\, is a continuous $K$-biframe for \,$H$\, with bounds \,$A,\,B$\, and \,$T \,\in\, \mathcal{B}\,(\,H\,)$\, be an invertible with \,$T\,K \,=\, K\,T$.\,Then \,$(\,T\,\mathcal{F},\, T\,\mathcal{G}\,)$\, is a continuous $K$-biframe for \,$H$.
\end{theorem}

\begin{proof}
For each \,$f \,\in\, H$, we have
\begin{align}
 \left\|\,K^{\,\ast}\,f \,\right\|^{\,2} &\,=\, \left\|\,\left (\,T^{\,-\, 1}\,\right )^{\,\ast}\,T^{\,\ast}\,K^{\,\ast}\,f\,\right\|^{\,2} \,\leq\, \left\|\,T^{\,-\, 1}\,\right\|^{\,2} \,\left\|\,T^{\,\ast}\,K^{\,\ast}\,f\,\right \|^{\,2}.\label{3.eq3.21}
\end{align}
Now, for each \,$f \,\in\, H$, we have
\begin{align*}
&\int\limits_{\Omega}\,\left<\,f,\, T\,F\,(\,w\,)\,\right>\,\left<\,T\,G\,(\,w\,),\, f\,\right>\,d\mu  \,=\, \int\limits_{\Omega}\,\left<\,T^{\,\ast}\,f,\, F\,(\,w\,)\,\right>\,\left<\,G\,(\,w\,),\, T^{\,\ast}\,f\,\right>\,d\mu \\
& \,\geq\, A\, \left\|\,K^{\,\ast}\,T^{\,\ast}\,f\,\right\|^{\,2} \,=\, A\,\left\|\,T^{\,\ast}\,K^{\,\ast}\,f \,\right\|^{\,2}\; \left[\;\text{since}\; T\,K \,=\, K\,T\,\right]\\
& \geq\, A\;\left \|\,T^{\,-\, 1}\,\right\|^{\,-2}\,\left\|\, K^{\,\ast}\,f\,\right\|^{\,2}\; \;[\;\text{using} \;(\ref{3.eq3.21})].
\end{align*}  
Following the proof of the Theorem \ref{3.ttm3.20}, it can be shown that \,$(\,T\,\mathcal{F},\, T\,\mathcal{G}\,)$\, satisfies upper frame condition.\,Hence, \,$(\,T\,\mathcal{F},\, T\,\mathcal{G}\,)$\, is a continuous $K$-biframe for \,$H$.
\end{proof}

\begin{corollary}
Let \,$(\,\mathcal{F},\, \mathcal{G}\,)$\, be a continuous $K$-biframe for \,$H$.\,Suppose \,$T \,\in\, \mathcal{B}(\,H\,)$\, such that \,$T\,T^{\,\ast} \,=\, I_{H}$ with \,$T\,K \,=\, K\,T$.\,Then \,$(\,T\,\mathcal{F},\, T\,\mathcal{G}\,)$\, is a continuous $K$-biframe for \,$H$.  
\end{corollary}

\begin{theorem}
Let \,$(\,\mathcal{F},\, \mathcal{G}\,)$\, be a continuous $K$-biframe for \,$H$\, with frame operator \,$S_{\mathcal{F},\,\mathcal{G}}$\, and \,$T$\, be a positive operator.\,Then \,$\left(\,\mathcal{F} \,+\, T\,\mathcal{F},\, \mathcal{G} \,+\, T\,\mathcal{G}\,\right)$\, is a continuous $K$-biframe for \,$H$.\,Furthermore, for any natural number \,$n$, \,$\left(\,\mathcal{F} \,+\, T^{\,n}\,\mathcal{F},\, \mathcal{G} \,+\, T^{\,n}\,\mathcal{G}\,\right)$\, is also continuous $K$-biframe for \,$H$. 
\end{theorem}

\begin{proof}
Let \,$(\,\mathcal{F},\, \mathcal{G}\,)$\, be a continuous $K$-biframe for \,$H$\, with frame operator \,$S_{\mathcal{F},\,\mathcal{G}}$.\,Then by Theorem \ref{3.thm3.8}, there exists \,$A \,>\, 0$\, such that \,$S_{\mathcal{F},\,\mathcal{G}} \,\geq\, A\,K\,K^{\,\ast}$.\,Now, for each \,$f \,\in\, H$, we have
\begin{align*}
&\int\limits_{\Omega}\,\left<\,f,\, F\,(\,w\,) \,+\, T\,F\,(\,w\,)\,\right>\,\left(\,G\,(\,w\,) \,+\, T\,G\,(\,w\,)\,\right)\,d\mu\\
&=\,\left(\,I_{H} \,+\, T\,\right)\int\limits_{\Omega}\,\left<\,f,\, \left(\,I_{H} \,+\, T\,\right)\,F\,(\,w\,)\,\right>\,G\,(\,w\,),\,d\mu \\
&=\,\left(\,I_{H} \,+\, T\,\right)\,S_{\mathcal{F},\,\mathcal{G}}\,\left(\,I_{H} \,+\, T\,\right)^{\,\ast}\,f.
\end{align*}
This shows that the frame operator for the frame \,$\left(\,\mathcal{F} \,+\, T\,\mathcal{F},\, \mathcal{G} \,+\, T\,\mathcal{G}\,\right)$\, is given by \,$\left(\,I_{H} \,+\, T\,\right)\,S_{\mathcal{F},\,\mathcal{G}}\,\left(\,I_{H} \,+\, T\,\right)^{\,\ast}$.\,Since \,$T$\, is positive, we get 
\[\left(\,I_{H} \,+\, T\,\right)\,S_{\mathcal{F},\,\mathcal{G}}\,\left(\,I_{H} \,+\, T\,\right)^{\,\ast} \,\geq\, S_{\mathcal{F},\,\mathcal{G}} \,\geq\, A\,K\,K^{\,\ast}.\]
Thus, \,$\left(\,\mathcal{F} \,+\, T\,\mathcal{F},\, \mathcal{G} \,+\, T\,\mathcal{G}\,\right)$\, is a continuous $K$-biframe for \,$H$. 

Furthermore, for any natural number \,$n$, it can be easily shown that the frame operator for \,$\left(\,\mathcal{F} \,+\, T^{\,n}\,\mathcal{F},\, \mathcal{G} \,+\, T^{\,n}\,\mathcal{G}\,\right)$\, is \,$\left(\,I_{H} \,+\, T^{\,n}\,\right)\,S_{\mathcal{F},\,\mathcal{G}}\,\left(\,I_{H} \,+\, T^{\,n}\,\right)^{\,\ast}\,\geq\, S_{\mathcal{F},\,\mathcal{G}} \,\geq\, A\,K\,K^{\,\ast}$.\,Hence, \,$\left(\,\mathcal{F} \,+\, T^{\,n}\,\mathcal{F},\, \mathcal{G} \,+\, T^{\,n}\,\mathcal{G}\,\right)$\, is a continuous $K$-biframe for \,$H$.  
\end{proof}

In the following Theorem, a characterization of continuous $K$-biframe is presented by using quotient of a bounded linear operator. 

\begin{theorem}
Let \,$K \,\in\, \mathcal{B}(\,H\,)$.\,Suppose \,$(\,\mathcal{F},\, \mathcal{G}\,)$\, is a continuous biframe Bessel mapping for \,$H$\, with frame operator \,$S_{\mathcal{F},\,\mathcal{G}}$.\,Then \,$(\,\mathcal{F},\, \mathcal{G}\,)$\, is a continuous $K$-biframe for \,$H$\, if and only if the quotient operator \,$\left[\,K^{\,\ast} \,/\, S^{1 \,/\, 2}_{\mathcal{F},\,\mathcal{G}}\,\right]$\, is bounded. 
\end{theorem}

\begin{proof}
First we suppose that \,$(\,\mathcal{F},\, \mathcal{G}\,)$\, is a continuous $K$-biframe for \,$H$\, with bounds \,$A$\, and \,$B$.\,Then for each \,$f \,\in\, H$, we have 
\[A\,\|\,K^{\,\ast}\,f\,\|^{\,2} \,\leq\, \int\limits_{\Omega}\,\left<\,f,\, F\,(\,w\,)\,\right>\,\left<\,G\,(\,w\,),\, f\,\right>\,d\mu  \,\leq\, B\,\|\,f \,\|^{\,2}.\]
Thus, for each \,$f \,\in\, H$, we have
\[A\,\|\,K^{\,\ast}\,f\,\|^{\,2} \,\leq\, \left<\,S_{\mathcal{F},\,\mathcal{G}}\,f,\, f\,\right> \,=\, \left\|\,S^{1 \,/\, 2}_{\mathcal{F},\,\mathcal{G}}\,f\,\right\|^{\,2}.\]
Now, it is easy to verify that the quotient operator \,$T \,:\, \mathcal{R}\left(\,S^{1 \,/\, 2}_{\mathcal{F},\,\mathcal{G}}\,\right) \,\to\, \mathcal{R}\,(\,K^{\,\ast}\,)$\, defined by \,$T\,\left(\,S^{1 \,/\, 2}_{\mathcal{F},\,\mathcal{G}}\,f\,\right) \,=\, K^{\,\ast}\,f\; \;\forall\; f \,\in\, H$\, is well-defined and bounded.

Conversely, suppose that the quotient operator \,$\left[\,K^{\,\ast} \,/\, S^{1 \,/\, 2}_{\mathcal{F},\,\mathcal{G}}\,\right]$\, is bounded.\,Then for each \,$f \,\in\, H$, there exists some \,$B \,>\, 0$\, such that
\begin{align*}
&\|\,K^{\,\ast}\,f\,\|^{\,2} \,\leq\, B\,\left\|\,S^{1 \,/\, 2}_{\mathcal{F},\,\mathcal{G}}\,f\,\right\|^{\,2} \,=\, B\,\left<\,S_{\mathcal{F},\,\mathcal{G}}\,f,\, f\,\right>\\
&\Rightarrow\, \|\,K^{\,\ast}\,f\,\|^{\,2} \,\leq\, B\,\int\limits_{\Omega}\,\left<\,f,\, F\,(\,w\,)\,\right>\,\left<\,G\,(\,w\,),\, f\,\right>\,d\mu.
\end{align*}
Thus, \,$(\,\mathcal{F},\, \mathcal{G}\,)$\, is a continuous $K$-biframe for \,$H$. 
\end{proof}

\begin{corollary}
Let \,$(\,\mathcal{F},\, \mathcal{G}\,)$\, be a continuous biframe Bessel mapping for \,$H$\, with frame operator \,$S_{\mathcal{F},\,\mathcal{G}}$.\,Then \,$(\,\mathcal{F},\, \mathcal{G}\,)$\, is a continuous biframe for \,$H$\, if and only if the frame operator  \,$S_{\mathcal{F},\,\mathcal{G}}$\, is invertible. 
\end{corollary}

Now, we establish that a quotient operator will be bounded if and only if a continuous $K$-biframe becomes continuous $T\,K$-biframe, for some \,$T \,\in\, \mathcal{B}\,(\,H\,)$.

\begin{theorem}
Let \,$(\,\mathcal{F},\, \mathcal{G}\,)$\, be a continuous $K$-biframe for \,$H$\, with frame operator \,$S_{\mathcal{F},\,\mathcal{G}}$\, and \,$T \,\in\, \mathcal{B}(\,H\,)$.\,Then the following are equivalent:
\begin{itemize}
\item[$(i)$]\,$\left(\,T\,\mathcal{F},\, T\,\mathcal{G}\,\right)$\, is a continuous $T\,K$-biframe for \,$H$.
\item[$(ii)$]\,$\left[\,\left(\,T\,K\,\right)^{\,\ast} \,/\, S^{1 \,/\, 2}_{\mathcal{F},\,\mathcal{G}}\,T^{\,\ast}\,\right]$\, is bounded. 
\item[$(iii)$]\,$\left[\,\left(\,T\,K\,\right)^{\,\ast} \,/\, \left(\,T\,S_{\mathcal{F},\,\mathcal{G}}\,T^{\,\ast}\,\right)^{1 \,/\, 2}\,\right]$\, is bounded. 
\end{itemize} 
\end{theorem}

\begin{proof}$(i) \,\Rightarrow\, (ii)$\, Suppose \,$\left(\,T\,\mathcal{F},\, T\,\mathcal{G}\,\right)$\, be a continuous $T\,K$-biframe for \,$H$\, with bounds \,$A$\, and \,$B$.\,Then for each \,$f \,\in\, H$, we have 
\[A\, \left\|\,\left(\,T\,K\,\right)^{\,\ast}\,f\,\right\|^{\,2} \,\leq\, \int\limits_{\Omega}\,\left<\,f,\, T\,F\,(\,w\,)\,\right>\,\left<\,T\,G\,(\,w\,),\, f\,\right>\,d\mu \,\leq\, B\, \|\,f\,\|^{\,2}.\]
Now, for each \,$f \,\in\, H$, we have
\begin{align}
\int\limits_{\Omega}\,\left<\,f,\, T\,F\,(\,w\,)\,\right>\,\left<\,T\,G\,(\,w\,),\, f\,\right>\,d\mu &\,=\, \int\limits_{\Omega}\,\left<\,T^{\,\ast}\,f,\, F\,(\,w\,)\,\right>\,\left<\,G\,(\,w\,),\, T^{\,\ast}\,f\,\right>\,d\mu \nonumber\\
& \,=\, \left<\,S_{\mathcal{F},\,\mathcal{G}}\,T^{\,\ast}\,f,\, T^{\,\ast}\,f\,\right>\label{eq3.6}.
\end{align}
Thus, for each \,$f \,\in\, H$, we have
\[A\, \left\|\,\left(\,T\,K\,\right)^{\,\ast}\,f\,\right\|^{\,2} \,\leq\, \left<\,S_{\mathcal{F},\,\mathcal{G}}\,T^{\,\ast}\,f,\, T^{\,\ast}\,f\,\right> \,=\, \left\|\,S^{1 \,/\, 2}_{\mathcal{F},\,\mathcal{G}}\,T^{\,\ast}\,f\,\right\|^{\,2}.\]
Now, we define a operator \,$Q \,:\, \mathcal{R}\left(\,S^{1 \,/\, 2}_{\mathcal{F},\,\mathcal{G}}\,T^{\,\ast}\,\,\right) \,\to\, \mathcal{R}\left(\,(\,T\,K\,)^{\,\ast}\,\right)$\, by \,$Q\,\left(\,S^{1 \,/\, 2}_{\mathcal{F},\,\mathcal{G}}\,T^{\,\ast}\,f\,\right) \,=\, (\,T\,K\,)^{\,\ast}\,f\; \;\forall\; f \,\in\, H$.\,Then it is easy verify that the quotient operator \,$Q$\, is well-defined and bounded. 

$(ii) \,\Rightarrow\, (iii)$\,\,It is obvious.

$(iii) \,\Rightarrow\, (i)$\, Suppose that the quotient operator \,$\left[\,\left(\,T\,K\,\right)^{\,\ast} \,/\, \left(\,T\,S_{\mathcal{F},\,\mathcal{G}}\,T^{\,\ast}\,\right)^{1 \,/\, 2}\,\right]$\, is bounded.\,Then for each \,$f \,\in\, H$, there exists \,$B \,>\, 0$\, such that 
\[\left\|\,\left(\,T\,K\,\right)^{\,\ast}\,f\,\right\|^{\,2} \,\leq\, B\, \left\|\,\left(\,T\,S_{\mathcal{F},\,\mathcal{G}}\,T^{\,\ast}\,\right)^{1 \,/\, 2}\,f\,\right\|^{\,2}.\]
Now, by (\ref{eq3.6}), for each \,$f \,\in\, H$, we have
\begin{align*}
&\int\limits_{\Omega}\,\left<\,f,\, T\,F\,(\,w\,)\,\right>\,\left<\,T\,G\,(\,w\,),\, f\,\right>\,d\mu \,=\, \left<\,S_{\mathcal{F},\,\mathcal{G}}\,T^{\,\ast}\,f,\, T^{\,\ast}\,f\,\right>\\
&=\,\left\|\,\left(\,T\,S_{\mathcal{F},\,\mathcal{G}}\,T^{\,\ast}\,\right)^{1 \,/\, 2}\,f\,\right\|^{\,2}\,\geq\, \dfrac{1}{B}\,\left\|\,\left(\,T\,K\,\right)^{\,\ast}\,f\,\right\|^{\,2}.
\end{align*}
Following the proof of the Theorem \ref{3.ttm3.20}, it can be shown that \,$(\,T\,\mathcal{F},\, T\,\mathcal{G}\,)$\, also satisfies upper frame condition.\,Hence, \,$\left(\,T\,\mathcal{F},\, T\,\mathcal{G}\,\right)$\, be a continuous $T\,K$-biframe for \,$H$.\,This completes the proof. 
\end{proof}

\begin{corollary}
Let \,$K \,\in\, \mathcal{B}(\,H\,)$\, and \,$(\,\mathcal{F},\, \mathcal{G}\,)$\, is a continuous biframe for \,$H$.\,Then following are equivalent:
\begin{itemize}
\item[$(i)$]\,$\left(\,K\,\mathcal{F},\, K\,\mathcal{G}\,\right)$\, be a continuous $K$-biframe for \,$H$.
\item[$(ii)$]\,$\left[\,K^{\,\ast} \,/\, S^{1 \,/\, 2}_{\mathcal{F},\,\mathcal{G}}\,\right]$\, is bounded.
\end{itemize}
\end{corollary}

\section{Continuous $K_{1} \,\otimes\, K_{2}$-biframe in $H_{1} \,\otimes\, H_{2}$}

In this section, we introduce the concept of continuous $K$-biframe in tensor product of Hilbert spaces \,$H_{1} \,\otimes\, H_{2}$\, and give a characterization.

\begin{definition}
Let \,$(\,X,\, \mu\,) \,=\, \left(\,X_{1} \,\times\, X_{2},\, \mu_{\,1} \,\otimes\, \mu_{\,2}\,\right)$\, be the product of measure spaces with \,$\sigma$-finite positive measures \,$\mu_{\,1},\, \mu_{\,2}$\, on \,$X_{1}$, \,$X_{2}$, respectively.\,Suppose \,$K_{1} \,\otimes\, K_{2}$\, be a bounded linear operator on \,$H_{1} \,\otimes\, H_{2}$.\,A pair mappings \,$(\,\mathcal{F},\, \mathcal{G}\,) \,=\, \,\left(\,F \,:\, X \,\to\, H_{1} \,\otimes\, H_{2},\, \,G \,:\, X \,\to\, H_{1} \,\otimes\, H_{2}\,\right)$\, is called a continuous $K_{1} \,\otimes\, K_{2}$-biframe for \,$H_{1} \,\otimes\, H_{2}$\, with respect to \,$(\,X,\, \mu\,)$\, if
\begin{itemize}
\item[$(i)$]$ F,\, G$\, is weakly-measurable, i.\,e., for all \,$f \,\otimes\, g \,\in\, H_{1} \,\otimes\, H_{2}$, \,$x \,=\, \left(\,x_{\,1},\, x_{\,2}\,\right) \,\mapsto\, \left<\,f \,\otimes\, g,\, F\,(\,x\,)\,\right>$\, and \,$\left(\,x_{\,1},\, x_{\,2}\,\right) \,\mapsto\, \left<\,f \,\otimes\, g,\, G\,(\,x\,)\,\right>$\, are measurable functions on \,$X$,
\item[$(ii)$]there exist constants \,$A,\,B \,>\, 0$\, such that
\begin{align}
&A\,\left\|\,\left(\,K_{1} \,\otimes\, K_{2}\,\right)^{\,\ast}\,(\,f \,\otimes\, g\,)\,\right\|^{\,2}\nonumber\\
& \,\leq\, \int\limits_{\,X}\,\left<\,f \,\otimes\, g,\, F\,(\,x\,)\,\right>\,\left<\,G\,(\,x\,),\, f \,\otimes\, g\,\right>\,d\mu\nonumber\\
&\leq\,B\,\left\|\,f \,\otimes\, g\,\right\|^{\,2},\label{4.eqq4.11} 
\end{align}
for all \,$f \,\otimes\, g \,\in\, H_{1} \,\otimes\, H_{2}$.\,The constants \,$A$\, and \,$B$\, are called frame bounds.\,If \,$A \,=\, B$, then the pair $(\,\mathcal{F},\, \mathcal{G}\,)$\, is called a tight continuous $K_{1} \,\otimes\, K_{2}$-biframe for \,$H_{1} \,\otimes\, H_{2}$. 
\end{itemize}   
\end{definition}

\begin{theorem}\label{th4.11}
The pair of mappings \,$\left(\,\mathcal{F},\, \mathcal{G}\,\right) \,=\,  \left(\,F_{1} \,\otimes\, F_{2},\, G_{1} \,\otimes\, G_{2}\,\right) \,=\, F,\, G \,:\, X \,\to\, H_{1} \,\otimes\, H_{2}$\, is a continuous $K_{1} \,\otimes\, K_{2}$-biframe for \,$H_{1} \,\otimes\, H_{2}$\, with respect to \,$(\,X,\, \mu\,)$\, if and only if \,$F_{1},\,G_{1} \,:\, X_{1} \,\to\, H_{1}$\, is a continuous $K_{1}$-biframe for \,$H_{1}$\, with respect to \,$\left(\,X_{1},\, \mu_{\,1}\,\right)$\, and \,$F_{2},\, G_{2} \,:\, X_{2} \,\to\, H_{2}$\, is a continuous $K_{2}$-biframe for \,$H_{2}$\, with respect to \,$\left(\,X_{2},\, \mu_{\,2}\,\right)$.  
\end{theorem}
 
\begin{proof}
Suppose that \,$\left(\,\mathcal{F},\, \mathcal{G}\,\right) \,=\,  \left(\,F_{1} \,\otimes\, F_{2},\, G_{1} \,\otimes\, G_{2}\,\right)$\, is a continuous $K_{1} \,\otimes\, K_{2}$-biframe for \,$H_{1} \,\otimes\, H_{2}$\, with respect to \,$(\,X,\, \mu\,)$\, having bounds \,$A$\, and \,$B$.\,Let \,$f \,\in\, H_{1} \,-\, \{\,\theta\,\}$\, and fix \,$g \,\in\, H_{2} \,-\, \{\,\theta\,\}$.\,Then \,$f \,\otimes\, g \,\in\, H_{1} \,\otimes\, H_{2} \,-\, \{\,\theta\,\otimes\,\theta\,\}$\, and by Fubini's theorem we have
\begin{align*}
&\int\limits_{\,X}\,\left<\,f \,\otimes\, g,\, F_{1}\,(\,x_{\,1}\,) \,\otimes\, F_{2}\,(\,x_{\,2}\,)\,\right>\,\left<\, G_{1}\,(\,x_{\,1}\,) \,\otimes\, G_{2}\,(\,x_{\,2}\,),\, f \,\otimes\, g\,\right>\,d\mu\\
&=\,\int\limits_{\,X_{1}}\,\left<\,f,\, F_{1}\,(\,x_{\,1}\,)\,\right>_{1}\,\left<\,G_{1}\,(\,x_{\,1}\,),\, f\,\right>_{1}\,d\mu_{\,1}\,\int\limits_{\,X_{2}}\,\left<\,g,\, F_{2}\,(\,x_{\,2}\,)\,\right>_{2}\,\left<\,G_{2}\,(\,x_{\,2}\,),\, g\,\right>_{2}\,d\mu_{\,2}.
\end{align*}
Therefore, for each \,$f \,\otimes\, g \,\in\, H_{1} \,\otimes\, H_{2}$, the inequality (\ref{4.eqq4.11}) can be written as
\begin{align*}
&A\,\left\|\,K^{\,\ast}_{1}\,f\,\right\|^{\,2}_{1}\,\left\|\,K^{\,\ast}_{2}\,g\,\right\|^{\,2}_{2}\\
&\leq\,\int\limits_{\,X_{1}}\,\left<\,f,\, F_{1}\,(\,x_{\,1}\,)\,\right>_{1}\,\left<\,G_{1}\,(\,x_{\,1}\,),\, f\,\right>_{1}\,d\mu_{\,1}\,\int\limits_{\,X_{2}}\,\left<\,g,\, F_{2}\,(\,x_{\,2}\,)\,\right>_{2}\,\left<\,G_{2}\,(\,x_{\,2}\,),\, g\,\right>_{2}\,d\mu_{\,2}\\
&\,\leq\,B\,\left\|\,f\,\right\|^{\,2}_{1}\,\left\|\,g\,\right\|^{\,2}_{2}.
\end{align*} 
Since \,$f$\, and \,$g$\, are non-zero and therefore 
\[\int\limits_{\,X_{1}}\,\left<\,f,\, F_{1}\,(\,x_{\,1}\,)\,\right>_{1}\,\left<\,G_{1}\,(\,x_{\,1}\,),\, f\,\right>_{1}\,d\mu_{\,1}\,,\, \int\limits_{\,X_{2}}\,\left<\,g,\, F_{2}\,(\,x_{\,2}\,)\,\right>_{2}\,\left<\,G_{2}\,(\,x_{\,2}\,),\, g\,\right>_{2}\,d\mu_{\,2}\] are non-zero.\,Thus from the above inequality we can write
\begin{align*}
&\dfrac{A\,\left\|\,K^{\,\ast}_{2}\,g\,\right\|^{\,2}_{2}}{\int\limits_{\,X_{2}}\,\left<\,g,\, F_{2}\,(\,x_{\,2}\,)\,\right>_{2}\,\left<\,G_{2}\,(\,x_{\,2}\,),\, g\,\right>_{2}\,d\mu_{\,2}}\,\left\|\,K^{\,\ast}_{1}\,f\,\right\|^{\,2}_{1}\\
&\leq\,\int\limits_{\,X_{1}}\,\left<\,f,\, F_{1}\,(\,x_{\,1}\,)\,\right>_{1}\,\left<\,G_{1}\,(\,x_{\,1}\,),\, f\,\right>_{1}\,d\mu_{\,1}\\
&\leq\,\dfrac{B\,\left\|\,g\,\right\|^{\,2}_{2}}{\int\limits_{\,X_{2}}\,\left<\,g,\, F_{2}\,(\,x_{\,2}\,)\,\right>_{2}\,\left<\,G_{2}\,(\,x_{\,2}\,),\, g\,\right>_{2}\,d\mu_{\,2}}\,\left\|\,f\,\right\|^{\,2}_{1}.
\end{align*}
Thus, for each \,$f \,\in\, H_{1} \,-\, \{\,\theta\,\}$, we have
\begin{align*}
A_{\,1}\,\left\|\,K^{\,\ast}_{1}\,f\,\right\|^{\,2}_{1} &\,\leq\, \int\limits_{\,X_{1}}\,\left<\,f,\, F_{1}\,(\,x_{\,1}\,)\,\right>_{1}\,\left<\,G_{1}\,(\,x_{\,1}\,),\, f\,\right>_{1}\,d\mu_{\,1} \leq\,B_{\,1}\,\left\|\,f\,\right\|^{\,2}_{1},
\end{align*}
where 
\[A_{\,1} \,=\, \inf_{g \,\in\, H_{\,2},\, \|\,g\,\|_{2} \,=\, 1}\,\left\{\,\dfrac{A\,\left\|\,K^{\,\ast}_{2}\,g\,\right\|^{\,2}_{2}}{\int\limits_{\,X_{2}}\,\left<\,g,\, F_{2}\,(\,x_{\,2}\,)\,\right>_{2}\,\left<\,G_{2}\,(\,x_{\,2}\,),\, g\,\right>_{2}\,d\mu_{\,2}}\,\right\},\]and
\[B_{\,1} \,=\, \sup_{g \,\in\, H_{\,2},\, \|\,g\,\|_{2} \,=\, 1}\,\left\{\,\dfrac{B\,\left\|\,g\,\right\|^{\,2}_{2}}{\int\limits_{\,X_{2}}\,\left<\,g,\, F_{2}\,(\,x_{\,2}\,)\,\right>_{2}\,\left<\,G_{2}\,(\,x_{\,2}\,),\, g\,\right>_{2}\,d\mu_{\,2}}\,\right\}.\] 
This shows that \,$\left(\,F_{1},\,G_{1}\,\right)$\, is a continuous $K_{1}$-biframe for \,$H_{1}$\, with respect to \,$\left(\,X_{1},\, \mu_{\,1}\,\right)$. Similarly, it can be shown that \,$\left(\,F_{2},\,G_{2}\,\right)$\, is a continuous $K_{2}$-biframe for \,$H_{2}$\, with respect to \,$\left(\,X_{2},\, \mu_{\,2}\,\right)$.

Conversely, suppose that \,$\left(\,F_{1},\,G_{1}\,\right)$\, is a continuous $K_{1}$-biframe for \,$H_{1}$\, with respect to \,$\left(\,X_{1},\, \mu_{\,1}\,\right)$\, having bounds \,$A,\, B$\, and \,$\left(\,F_{2},\,G_{2}\,\right)$\, is a continuous $K_{2}$-biframe for \,$H_{2}$\, with respect to \,$\left(\,X_{2},\, \mu_{\,2}\,\right)$\, having bounds \,$C,\,D$.\,By the assumption it is easy to very that \,$\left(\,\mathcal{F},\, \mathcal{G}\,\right) \,=\,  \left(\,F_{1} \,\otimes\, F_{2},\, G_{1} \,\otimes\, G_{2}\,\right)$ is weakly measurable on \,$H_{1} \,\otimes\, H_{2}$\, with respect to \,$(\,X,\, \mu\,)$.\,Now, for each \,$f \,\in\, H_{1} \,-\, \{\,\theta\,\}$\, and \,$g \,\in\, H_{2} \,-\, \{\,\theta\,\}$, we have
\begin{align*}
&A\,\left\|\,K^{\,\ast}_{1}\,f\,\right\|^{\,2}_{1} \,\leq\, \int\limits_{\,X_{1}}\,\left<\,f,\, F_{1}\,(\,x_{\,1}\,)\,\right>_{1}\,\left<\,G_{1}\,(\,x_{\,1}\,),\, f\,\right>_{1}\,d\mu_{\,1} \leq\,B\,\left\|\,f\,\right\|^{\,2}_{1}\,,\\
&C\,\left\|\,K^{\,\ast}_{2}\,g\,\right\|^{\,2}_{2} \,\leq\, \int\limits_{\,X_{2}}\,\left<\,g,\, F_{2}\,(\,x_{\,2}\,)\,\right>_{2}\,\left<\,G_{2}\,(\,x_{\,2}\,),\, g\,\right>_{1}\,d\mu_{\,2} \leq\,B\,\left\|\,g\,\right\|^{\,2}_{2}.
\end{align*} 
Multiplying the above two inequalities and using Fubini's theorem we get
\begin{align*}
&A\,C\,\left\|\,\left(\,K_{1} \,\otimes\, K_{2}\,\right)^{\,\ast}\,(\,f \,\otimes\, g\,)\,\right\|^{\,2}\\
& \,\leq\, \int\limits_{\,X}\,\left<\,f \,\otimes\, g,\, F_{1}\,(\,x_{\,1}\,) \,\otimes\, F_{2}\,(\,x_{\,2}\,)\,\right>\,\left<\, G_{1}\,(\,x_{\,1}\,) \,\otimes\, G_{2}\,(\,x_{\,2}\,),\, f \,\otimes\, g\,\right>\,d\mu \\
&\leq\,B\,D\,\left\|\,f \,\otimes\, g\,\right\|^{\,2}, 
\end{align*} 
for all \,$f \,\otimes\, g \,\in\, H_{1} \,\otimes\, H_{2}$.\,Thus, for each \,$f \,\otimes\, g \,\in\, H_{1} \,\otimes\, H_{2}$, we have
\begin{align*}
A\,C\,\left\|\,\left(\,K_{1} \,\otimes\, K_{2}\,\right)^{\,\ast}\,(\,f \,\otimes\, g\,)\,\right\|^{\,2} &\,\leq\, \int\limits_{\,X}\,\left<\,f \,\otimes\, g,\, F\,(\,x\,)\,\right>\,\left<\,G\,(\,x\,),\, f \,\otimes\, g\,\right>\,d\mu \\
&\leq\,B\,D\,\left\|\,f \,\otimes\, g\,\right\|^{\,2},
\end{align*}
Hence, \,$\left(\,\mathcal{F},\, \mathcal{G}\,\right) \,=\,  \left(\,F_{1} \,\otimes\, F_{2},\, G_{1} \,\otimes\, G_{2}\,\right)$\, is a continuous $K_{1} \,\otimes\, K_{2}$-biframe for \,$H_{1} \,\otimes\, H_{2}$\, with respect to \,$(\,X,\, \mu\,)$\, having bounds \,$A\,C$\, and \,$B\,D$.\,This completes the proof.  
\end{proof} 

\begin{example}
Let \,$\left\{\,e_{\,i}\,\right\}_{i \,=\, 1}^{\,\infty}$\, be an orthonormal basis for \,$H_{1}$\, and \,$\left(\,X_{1},\, \mu_{\,1}\,\right)$\, be a measure space with \,$\mu_{\,1}$\, is \,$\sigma$-finite.\,Then we can write \,$X_{1} \,=\, \bigcup_{i \,=\, 1}^{\,\infty}\,\Omega_{i}$, where \,$\left\{\,\Omega_{\,i}\,\right\}_{i \,=\, 1}^{\,\infty}$\, is a sequence of disjoint measurable subsets of \,$X_{1}$\, with \,$\mu\left(\,\Omega_{i}\,\right) \,<\, \infty$.\,Suppose
\begin{align*}
&\left\{\,f_{\,i}\,\right\}_{i \,=\, 1}^{\,\infty} \,=\, \left\{\,3\,e_{\,1},\, 2\,e_{\,2},\, 2\,e_{\,3},\, \cdots\,\cdots\,\right\}\,,\\
&\left\{\,g_{\,i}\,\right\}_{i \,=\, 1}^{\,\infty} \,=\, \left\{\,0,\, 2\,e_{\,1},\, 0,\, e_{\,2},\, 0,\, e_{\,3},\, \cdots\,\cdots\,\right\}.
\end{align*}
Define \,$K_{1} \,:\, H_{1} \,\to\, H_{1}$\, by \,$K_{1}\,f \,=\, \sum\limits_{i \,=\, 1}^{\,m}\,\left<\,f,\, e_{\,i}\,\right>\,e_{\,i}$, \,$f \,\in\, H_{1}$, where \,$m$\, is a fixed positive integer.\,For each \,$x_{1} \,\in\, X_{1}$, we define the mappings \,$F_{1} \,:\, X_{1} \,\to\, H_{1}$\, by \,$F_{1}\,(\,x_{1}\,) \,=\, \dfrac{1}{\sqrt{\mu\left(\,\Omega_{i}\,\right)}}\,f_{\,i}$\, and \,$G_{1} \,:\, X_{1} \,\to\, H_{1}$\, by \,$G_{1}\,(\,x_{1}\,) \,=\, \dfrac{1}{\sqrt{\mu\left(\,\Omega_{i}\,\right)}}\,g_{\,i}$.
Then by Example \ref{3exm3.11}, \,$\left(\,F_{1},\, G_{1}\,\right)$\, is a continuous $K_{1}$-biframe for \,$H_{1}$\, with bounds \,$1$\, and \,$2$.

On the other hand, let \,$\left\{\,e^{\,\prime}_{\,j}\,\right\}_{j \,=\, 1}^{\,\infty}$\, be an orthonormal basis for \,$H_{2}$\, and \,$\left(\,X_{2},\, \mu_{\,2}\,\right)$\, be a measure space with \,$\mu_{\,2}$\, is \,$\sigma$-finite.\,Then \,$X_{2} \,=\, \bigcup_{i \,=\, 1}^{\,\infty}\,\Omega^{\,\prime}_{j}$, and \,$\left\{\,\Omega^{\,\prime}_{\,j}\,\right\}_{j \,=\, 1}^{\,\infty}$\, is a sequence of disjoint measurable subsets of \,$X_{2}$\, with \,$\mu\left(\,\Omega^{\,\prime}_{j}\,\right) \,<\, \infty$.\,Suppose 
\begin{align*}
&\left\{\,f^{\,\prime}_{\,j}\,\right\}_{i \,=\, 1}^{\,\infty} \,=\, \left\{\,5\,e_{\,1},\, 3\,e_{\,2},\, 2\,e_{\,3},\, 2\,e_{\,4},\, \cdots\,\cdots\,\right\}\,,\\
&\left\{\,g^{\,\prime}_{\,j}\,\right\}_{i \,=\, 1}^{\,\infty} \,=\, \left\{\,0,\, 0,\, 3\,e_{\,1},\, 0,\, 2\,e_{\,2},\, 2\,e_{\,3},\, \cdots\,\cdots\,\right\}.
\end{align*}
Define \,$K_{2} \,:\, H_{2} \,\to\, H_{2}$\, by \,$K_{2}\,g \,=\, \sum\limits_{j \,=\, 1}^{\,\infty}\,\left<\,g,\, e_{\,2\,j}\,\right>\,e_{\,2\,j}$, \,$g \,\in\, H_{2}$.\,Then for \,$g \,\in\, H_{2}$, we get
\[\left\|\,K^{\,\ast}_{2}\,g\,\right\|^{\,2}_{2} \,=\, \sum\limits_{j \,=\, 1}^{\,\infty}\,\left|\,\left<\,g,\, e_{\,2\,j}\,\right>\,\right|^{\,2} \,\leq\, \|\,g\,\|_{2}^{\,2}.\]
Now, we define \,$F_{\,2} \,:\, X_{2} \,\to\, H_{2}$\, by \,$F_{\,2}\,(\,x_{\,2}\,) \,=\, \dfrac{1}{\sqrt{\mu\left(\,\Omega^{\,\prime}_{j}\,\right)}}\,f^{\,\prime}_{\,j}$\, and  \,$G_{\,2} \,:\, X_{2} \,\to\, H_{2}$\, by \,$G_{\,2}\,(\,x_{\,2}\,) \,=\, \dfrac{1}{\sqrt{\mu\left(\,\Omega^{\,\prime}_{j}\,\right)}}\,g^{\,\prime}_{\,j}$.\,Now, for \,$g \,\in\, H_{2}$, we have
\begin{align*}
&\int\limits_{\,X_{2}}\,\left<\,g,\, F_{2}\,(\,x_{2}\,)\,\right>\,\left<\,G_{2}\,(\,x_{2}\,),\, g\,\right>\,d\mu_{2} \,=\,\sum\limits_{j \,=\, 1}^{\infty}\,\int\limits_{\,\Omega^{\,\prime}_{j}}\,\left <\,g,\, f^{\,\prime}_{\,j}\,\right >\,\left <\,g^{\,\prime}_{\,j},\, g\,\right >\,d\mu_{2}\\
&=\,\left <\,g,\, e_{\,1}\,\right >\,\left <\,e_{\,1},\, g\,\right > \,+\, 2\,\left <\,g,\, e_{\,1}\,\right >\,\left <\,e_{\,1},\, g\,\right > \,+\, 2\,\left <\,g,\, e_{\,2}\,\right >\,\left <\,e_{\,2},\, g\,\right >+\cdots\\
&=\,\left|\,\left <\,g,\, e_{\,1}\,\right >\,\right|^{\,2} \,+\, 2\,\left|\,\left <\,g,\, e_{\,1}\,\right >\,\right|^{\,2} \,+\, 2\,\left|\,\left <\,g,\, e_{\,2}\,\right >\,\right|^{\,2} \,+\, \cdots\\
& \,=\, \left|\,\left <\,g,\, e_{\,1}\,\right >\,\right|^{\,2} \,+\, 2\,\|\,g\,\|^{\,2}. 
\end{align*} 
Thus, \,$\left(\,F_{2},\, G_{2}\,\right)$\, is a continuous $K_{2}$-biframe for \,$H$\, with bounds \,$1$\, and \,$3$.\,Thus, by the Theorem \ref{th4.11}, \,$\left(\,\mathcal{F},\, \mathcal{G}\,\right) \,=\,  \left(\,F_{1} \,\otimes\, F_{2},\, G_{1} \,\otimes\, G_{2}\,\right)$\, is a continuous $K_{1} \,\otimes\, K_{2}$-biframe for \,$H_{1} \,\otimes\, H_{2}$\, with respect to \,$(\,X,\, \mu\,)$\, having bounds \,$1$\, and \,$6$.
\end{example}

\end{document}